\documentclass[onecolumn,a4paper,12pt,oneside]{article} 
\usepackage[top=2.5cm,left=2.5cm,right=2.5cm,bottom=3cm]{geometry}    
\usepackage{draftwatermark} 
\SetWatermarkText{\bf Preprint}
\SetWatermarkAngle{45}
\SetWatermarkScale{3}
\SetWatermarkLightness{0.85}     
\usepackage[parfill]{parskip}    
\usepackage{graphicx}
\usepackage{amssymb}
\usepackage{amsmath}
\usepackage{amsfonts}
\usepackage{amsthm}
\usepackage{epstopdf}
\usepackage[usenames,dvipsnames,svgnames,table]{xcolor}
\usepackage{tikz}
\usepackage{array}
\usepackage[T1]{fontenc}
\usepackage[utf8]{inputenc}
\usepackage{booktabs} 
\usepackage{array} 
\usepackage{paralist} 
\usepackage{listings} 
\usepackage{verbatim} 
\usepackage{subfig} 
\usepackage[colorlinks=true]{hyperref}
\hypersetup{colorlinks=true, linktocpage=true, linkcolor=Blue, citecolor=Green,
hypertexnames=true, urlcolor=MidnightBlue,pdftex}
\usepackage[english]{babel}
\usepackage{wrapfig}
\usepackage{multirow}
\usepackage{quoting}
\usepackage{rotating}
\quotingsetup{font=normalsize}
\theoremstyle{plain}

\usepackage{bm}
\usepackage{abstract}
\usepackage{float}
\usepackage{textcomp} 	
\usepackage{tabularx} 
\usepackage{caption}
\usepackage{authblk} 
\usepackage{eso-pic} 
\usepackage{lineno} 
\providecommand{\keywords}[1]{\textbf{{Key words: }} #1} 

\newcommand{\be}{\begin{equation}}
\newcommand{\ee}{\end{equation}}
\newcommand{\bs}{\begin{split}}
\newcommand{\es}{\end{split}}

\newcommand{\beps}{\boldsymbol{\varepsilon}}
\newcommand{\balpha}{\boldsymbol{\alpha}}

\newcommand{\bk}{\boldsymbol{\kappa}}
\newcommand{\bgamma}{\boldsymbol{\gamma}}

\renewcommand{\Phi}{\varPhi}
\newcommand{\gr}{\mathbf}

\newcommand{\tr}{\mathrm{tr}}

\renewcommand{\Theta}{\varTheta}
\renewcommand{\Psi}{\varPsi}
\renewcommand{\Sigma}{\varSigma}
\newcommand{\A}{\mathbb{A}}
\newcommand{\Ac}{\mathcal{A}}
\newcommand{\B}{\mathbb{B}}
\newcommand{\Bc}{\mathcal{B}}
\newcommand{\C}{\mathbb{C}}
\newcommand{\D}{\mathbb{D}}
\newcommand{\Dc}{\mathcal{D}}
\newcommand{\Q}{\mathbb{Q}}

\renewcommand{\Delta}{\varDelta}
\renewcommand{\phi}{\varphi}
\renewcommand{\psi}{\varPsi}

\newcommand{\N}{\mathbf{N}}
\newcommand{\M}{\mathbf{M}}
\newcommand{\U}{\mathbf{U}}
\newcommand{\V}{\mathbf{V}}
\newcommand{\W}{\mathbf{W}}

\renewcommand{\u}{\mathbf{u}}
\renewcommand{\v}{\mathbf{v}}
\newcommand{\w}{\mathbf{w}}

\title{\LARGE{\bf On the thermoelastic coupling of anisotropic laminates}}

\author{Paolo Vannucci\\ 
\begin{small}LMV - Laboratoire de Mathématiques de Versailles, UMR8100.\\Université de Versailles et Saint Quentin - \href{mailto:paolo.vannucci@uvsq.fr}{paolo.vannucci@uvsq.fr}\end{small}}

\begin{document}

\maketitle

PREPRINT - Final version: Archive of Applied Mechanics, v. 94: 1121-1149, 2024.\bigskip\\

\hrule
\begin{abstract}
The analysis of the mathematical and mechanical properties of thermoelastic coupling tensors in anisotropic laminates is the topic of this paper. Some theoretical results concerning the compliance tensors are shown and their mechanical consequences analyzed. Moreover, the case of thermally stable laminates, important for practical applications, is also considered. The study is carried out in the framework of the polar method, a mathematical formalism particularly well suited for the analysis of planar anisotropic problems, introduced by Prof. G. Verchery in 1979.

\end{abstract}
\medskip
\hrule
\bigskip



\section{Introduction}
\label{sec1}
Anisotropic laminates, obtained superposing layers of composite materials, can show a coupling between the extension and the bending behaviors. This coupling is due to the stacking sequence, if the layers are all identical, and also to the difference of the layers, in case of hybrid laminates, i.e. composed by layers of different materials. Normally, laminates are designed to be uncoupled but coupling can be an interesting device to obtain some peculiar mechanical responses of interest in some applications. The effects of coupling when the laminate is acted upon by mechanical actions has been the subject of a recent paper, \cite{vannucci23a}, where, however, the consequences of thermal loads on a coupled laminate were not considered: this is the topic of this paper. In particular, we analyze the mathematical properties of the tensors describing the thermoelastic coupling and its relations with the elastic couplings,   for both the stiffness and compliance tensors. We  also analyze the mechanical consequences of thermoelastic coupling in some cases of interest. The polar method, a mathematical technique introduced by Professor G. Verchery, \cite{verchery79,meccanica05,vannucci_libro} is used to represent all the tensors. In this method, any planar tensor, and its Cartesian components,  is represented using exclusively invariants and angles, which is a true advantage with anisotropic problems. In addition, with the polar method all the material symmetries are represented by special values of some tensor invariants, another important advantage for the mathematical statement of design problems in anisotropic materials and structures. 

In the first part of the paper, the governing equations are introduced and the essential points of the polar method, of interest for the paper, are recalled. Then, the algebraic aspects of the thermoelastic coupling tensors, in stiffness and in compliance, are analyzed and finally the mechanical response of a set of coupled laminates is studied.

\section{The governing equations}
In the classical theory of anisotropic laminated plates, \cite{jones,vannucci_libro}, the constitutive law is 
 \begin{equation}
 \label{eq:fundlawtherm}
 \left\{\begin{array}{c}\gr{N} \\\hline \gr{M}\end{array}\right\}=
 \left[\begin{array}{c|c}h\A & \frac{h^2}{2}\B \\\hline \frac{h^2}{2}\B & \frac{h^3}{12}\D\end{array}\right] \left\{\begin{array}{c}\beps \\\hline \bk\end{array}\right\}
 -t \left\{\begin{array}{c}h\gr{U} \\\hline \frac{h^2}{2}\gr{V}\end{array}\right\}-\nabla t \left\{\begin{array}{c}\frac{h^2}{2}\gr{V} \\\hline \frac{h^3}{12}\gr{W}\end{array}\right\},
 \end{equation}
with  $\N$ and $\M$  respectively the extension and bending internal actions tensors, $\beps$  the  tensor of the in-plane deformation, $\bk$ the curvature tensor, $h$ the laminate's thickness, while $t$ is the temperature variation with respect to the temperature corresponding to a no strain state and $\nabla t$ is the through-the-thickness gradient of temperature. The tensors describing the elastic behavior of the laminate are   $\A,\B$ and $\D$, respectively the extension, coupling and bending stiffness tensors, homogeneous to an elasticity tensor and having the same classes of symmetry, i.e. the minor,
\be
\A_{ijkl}=\A_{jikl}=\A_{ijlk},
\ee
and major symmetries (of the indexes):
\be
\A_{ijkl}=\A_{klij}.
\ee
The major symmetries are actually the condition for a fourth-rank tensor to be defined as symmetric: $\A=\A^\top$, cf. \cite{vannucci_alg}. Of course,  the same can be said also for $\B$ and $\D$. 

The stiffness tensors describing the thermoelastic response are  $\U,\V$ and $\W$, that are the analogous of $\A,\B$ and $\D$ respectively, but for a thermal action: $\U$ determines the in-plane forces due to a through-the-thickness uniform change of temperature $t$, $\W$ the couples linked to a temperature gradient  $\nabla t$ and $\V$ is the thermoelastic coupling stiffness  tensor: it determines the in-plane forces due to a temperature gradient $\nabla t$ and the couples  caused by a uniform change of temperature $t$. All of them are  second-rank symmetric tensors.

The general expressions of $\A,\B,\D$ are 
 \begin{equation}
 \label{eq:compABD}
 \begin{split}
 &\A=\frac{1}{h}\sum_{k=1}^n\int_{z_{k-1}}^{z_k}\Q(\delta_k)dx_3=\frac{1}{h}\sum_{k=1}^n({z_k}-{z_{k-1}})\Q(\delta_k),\\
 &\B=\frac{2}{h^2}\sum_{k=1}^n\int_{z_{k-1}}^{z_k}x_3\Q(\delta_k)dx_3=\frac{1}{h^2}\sum_{k=1}^n({z_k^2}-{z_{k-1}^2})\Q(\delta_k),\\
& \D=\frac{12}{h^3}\sum_{k=1}^n\int_{z_{k-1}}^{z_k}x_3^2\Q(\delta_k)dx_3=\frac{4}{h^3}\sum_{k=1}^n({z_k^3}-{z_{k-1}^3})\Q(\delta_k),
 \end{split}
 \end{equation}
with $\mathbb{Q}(\delta_k)$ the reduced stiffness tensor of the $k-$th layer, at the orientation $\delta_k$, $n$ the number of plies  and $z_k,z_{k-1}$ the position of the upper and lower surfaces of the same layer on the thickness of the plate. 
Similarly,
\be
\label{eq:tensorthemoelast}
\begin{split}
&\gr{U}=\frac{1}{h}\sum_{k=1}^n\int_{z_{k-1}}^{z_k}\bgamma_k(\delta_k)dx_3=\frac{1}{h}\sum_{k=1}^n(z_k-z_{k-1})\bgamma_k(\delta_k),\\
&\gr{V}=\frac{2}{h^2}\sum_{k=1}^n\int_{z_{k-1}}^{z_k}x_3\bgamma_k(\delta_k)dx_3=\frac{1}{h^2}\sum_{k=1}^n(z_k^2-z_{k-1}^2)\bgamma_k(\delta_k),\\
&\gr{W}=\frac{12}{h^3}\sum_{k=1}^n\int_{z_{k-1}}^{z_k}x_3^2\bgamma_k(\delta_k)dx_3=\frac{4}{h^3}\sum_{k=1}^n(z_k^3-z_{k-1}^3)\bgamma_k(\delta_k),
\end{split}
\ee
where the tensor $\bgamma(\delta_k)$ is given by
 \be
\label{eq:tensorgammatermico}
\bgamma(\delta_k)=\mathbb{Q}(\delta_k)\balpha(\delta_k), \ \ k=1,...,n,
\ee
with $\balpha(\delta_k)$ the tensor of the thermal expansion coefficients of $k-$th layer, at the orientation $\delta_k$.

The converse of the constitutive equation (\ref{eq:fundlawtherm}) is, cf. \cite{vannucci_libro},
 \begin{equation}
  \label{eq:fundlawinversetherm}
 \left\{\begin{array}{c}{\beps} \\\hline \bk\end{array}\right\}=
 \left[\begin{array}{c|c}\frac{1}{h}\mathcal{A} & \frac{2}{h^2}\mathcal{B} \\\hline \frac{2}{h^2}\mathcal{B}^\top & \frac{12}{h^3}\mathcal{D}\end{array}\right]
 \left\{\begin{array}{c}\gr{N} \\\hline \gr{M}\end{array}\right\}+t \left\{\begin{array}{c} \gr{u} \\\hline \gr{v}_1\end{array}\right\}+\nabla t \left\{\begin{array}{c}\gr{v}_2 \\\hline \gr{w}\end{array}\right\},
 \end{equation}
with $\Ac,\Bc$ and $\Dc$  the  compliance corresponding of $\A,\B$ and $\D$, given in the order by
  \begin{equation}
 \label{eq:inversetensorsABD}
 \begin{split}
 &\mathcal{A}=(\A-3\B\D^{-1}\B)^{-1},\\ 
 &\mathcal{B}=-3\mathcal{A}\B\D^{-1}=(-3\mathcal{D}\B\A^{-1})^\top=-3\A^{-1}\B\mathcal{D},\\
  &\mathcal{D}=(\D-3\B\A^{-1}\B)^{-1},
\end{split}
 \end{equation}
 and
\be
\label{eq:tenstherminverses}
\begin{split}
&\gr{u}=\mathcal{A}\gr{U}+\mathcal{B}\gr{V}=(\A-3\B\D^{-1}\B)^{-1}(\gr{U}-3\B\D^{-1}\gr{V}),\\
&\gr{v}_1=\frac{2}{h}(\mathcal{B}^\top\gr{U}+3\mathcal{D}\gr{V})=\frac{6}{h}(\D-3\B\A^{-1}\B)^{-1}(\gr{V}-\B\A^{-1}\gr{U}),\\
&\gr{v}_2=\frac{h}{6}(3\mathcal{A}\gr{V}+\mathcal{B}\gr{W})=\frac{h}{2}(\A-3\B\D^{-1}\B)^{-1}(\gr{V}-\B\D^{-1}\gr{W}),\\
&\gr{w}=\mathcal{B}^\top\gr{V}+\mathcal{D}\gr{W}=(\D-3\B\A^{-1}\B)^{-1}(\gr{W}-3\B\A^{-1}\gr{V}),
\end{split}
\ee
the  compliance thermoelastic tensors. Contrarily to what could seem, the compliance thermoelastic tensors are only three because 
\be
\label{eq:legamev1v2}
\gr{v}_2=\frac{h^2}{12}\mathcal{A}\left[\mathcal{D}^{-1}\gr{v}_1+\frac{6}{h}\B(\A^{-1}\gr{U}-\D^{-1}\gr{W})\right].
\ee
In particular:
\begin{itemize}
\item $\gr{u}$ is the tensor of the coefficients of thermal expansion of the laminate in case of a uniform temperature change $t$; its SI units are $^\circ$C$^{-1}$;
\item $\gr{v}_1$ is the tensor of the coefficients of thermal expansion of the laminate due to a gradient of temperature $\nabla t$; its SI units are $(m\ ^\circ$C$)^{-1}$;
\item $\gr{v}_2$ is the tensor of the coefficients of thermal curvature of the laminate for a uniform change of temperature $t$; its SI units are $m\ ^\circ$C$^{-1}$;
\item $ \gr{w}$ is the tensor of the coefficients of thermal curvature of the laminate caused by a gradient of temperature $\nabla t$; its SI units are $^\circ$C$^{-1}$.
\end{itemize}

The symmetries of elasticity and of $\balpha$ ensure the symmetry of $\U,\V$ and $\W$, which in turn, through the minor symmetries of $\Ac,\Bc,\Dc$, gives the symmetries of $\u,\v_1,\v_2,\w$.

We adopt the Kelvin's notation, \cite{kelvin,kelvin1}, for  2nd-rank tensors, e.g.:
\be
\N=\left\{
\begin{array}{c}
\N_1=\N_{11}\\
\N_2=\N_{22}\\
\N_6=\sqrt{2}\N_{12}
\end{array}
\right\},\ \
\ee
and similarly for $\M,\beps,\bk,\U,\V,\W,\u,\v_1,\v_2$ and $\w$, as well as for 4th-rank tensors, like
\be
\label{eq:kelvinmatrix}
\A\hspace{-1mm}=\hspace{-1mm}\left[
\begin{array}{ccc}
\A_{11}\hspace{-1mm}=\hspace{-1mm}\A_{1111}&\hspace{-1mm}\A_{12}\hspace{-1mm}=\hspace{-1mm}\A_{1122}&\hspace{-1mm}\A_{16}\hspace{-1mm}=\hspace{-1mm}\sqrt{2}\A_{1112}\\
\A_{12}\hspace{-1mm}=\hspace{-1mm}\A_{1122}&\hspace{-1mm}\A_{22}\hspace{-1mm}=\hspace{-1mm}\A_{2222}&\hspace{-1mm}\A_{26}\hspace{-1mm}=\hspace{-1mm}\sqrt{2}\A_{2212}\\
\A_{16}\hspace{-1mm}=\hspace{-1mm}\sqrt{2}\A_{1112}&\hspace{-1mm}\A_{26}\hspace{-1mm}=\hspace{-1mm}\sqrt{2}\A_{2212}&\hspace{-1mm}\A_{66}\hspace{-1mm}=\hspace{-1mm}2\A_{1212}
\end{array}
\right],
\ee
and the same for $\B,\D,\Ac,\Bc,\Dc$.

\section{Some basic elements of the polar method for laminates}
\begin{figure}[h]
\center
\includegraphics[width=.5\textwidth]{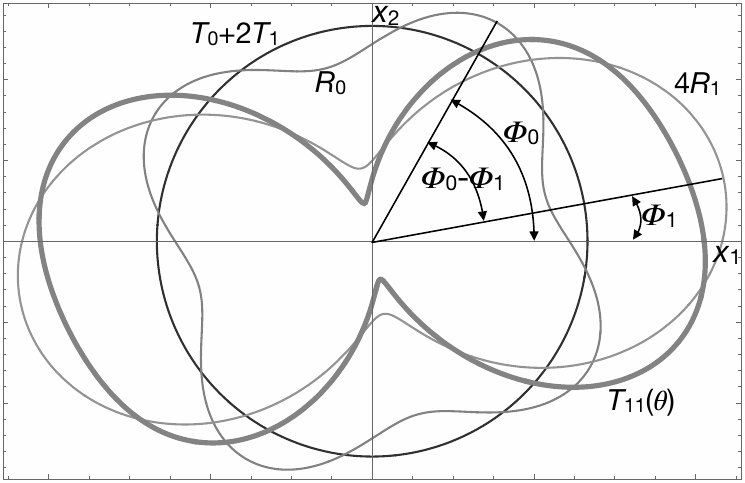}
\caption{The anisotropic phases of the component $\mathbb{T}_{11}(\theta)$ for a generic material.}
\label{fig:1}
\end{figure}

By the polar method, the Cartesian components in the Kelvin notation at a direction $\theta$ of a plane elastic tensor $\mathbb{T}$ are expressed as
\begin{equation}
\label{eq:mohr4}
{\begin{split}
&{\mathbb{T}_{11}{(}\theta{)}{=}{T}_{0}{+}{2}{T}_{1}{+}{R}_{0}\cos{4}\left({{{\varPhi}}_{0}{-}\theta}\right){+}{4}{R}_{1}\cos{2}\left({{{\varPhi}}_{1}{-}\theta}\right)},\\
&{\mathbb{T}_{12}{(}\theta{)}{=}{-}{T}_{0}{+}{2}{T}_{1}{-}{R}_{0}\cos{4}\left({{{\varPhi}}_{0}{-}\theta}\right)},\\
&\mathbb{T}_{16}(\theta)=\sqrt{2}\left[R_0\sin4\left(\varPhi_0-\theta\right)+2R_1\sin2\left(\varPhi_1-\theta\right)\right],\\
&{\mathbb{T}_{22}{(}\theta{)}{=}{T}_{0}{+}{2}{T}_{1}{+}{R}_{0}\cos{4}\left({{{\varPhi}}_{0}{-}\theta}\right){-}{4}{R}_{1}\cos{2}\left({{{\varPhi}}_{1}{-}\theta}\right)},\\
&{\mathbb{T}_{26}{(}\theta{)}{=}\sqrt{2}\left[{-}{R}_{0}\sin{4}\left({{{\varPhi}}_{0}{-}\theta}\right){+}{2}{R}_{1}\sin{2}\left({{{\varPhi}}_{1}{-}\theta}\right)\right]},\\
&{\mathbb{T}_{66}{(}\theta{)}{=}2\left[{T}_{0}{-}{R}_{0}\cos{4}\left({{{\varPhi}}_{0}{-}\theta}\right)\right]}.
\end{split}}
\end{equation}
The moduli $T_0,T_1,R_0,R_1$ as well as the difference of the angles $\Phi_0-\Phi_1$ are tensor invariants.  Moreover, to choose the value of one of the two polar angles corresponds to fix a frame; usually $\Phi_1=0$. The elastic symmetries are determined by the following conditions on the invariants:
\begin{itemize}
\item ordinary orthotropy: $\Phi_0-\Phi_1=k\dfrac{\pi}{4}, \ k\in\{0,1\}$;
\item $R_0$-orthotropy: $R_0=0$, \cite{vannucci02joe};
\item square symmetry: ${R_1=0}$;
\item isotropy: $R_0=R_1=0$.
\end{itemize}
So we see that $T_0$ and $T_1$ are the {\it isotropy invariants}, while $R_0,R_1$ and $\Phi_0-\Phi_1$ are the {\it anisotropy invariants}. In practice, in the polar method plane anisotropic elasticity is decomposed into an {\it isotropic phase}, depending on the two invariants $T_0$ and $T_1$, plus two {\it anisotropic phases}, shifted of the invariant angle $\Phi_0-\Phi_1$ and whose amplitudes are proportional to $R_0$ and $R_1$ respectively; the polar curves in Fig. \ref{fig:1} represent these elastic phases for a generic material.

The above relations are valid for any tensor of the elastic type, hence  for  $\A,\B,\D,\Ac,\Dc$ too, {\it but not for} $\Bc$, because it is not symmetric. 
In the following, a superscript $A,B$ or $D$ denotes a polar quantity of $\A,\B$ or $\D$ respectively; capital letters are reserved to stiffness tensors and lower-case letters to  compliance tensors (e.g., $T_0^A,\Phi_1^A$ etc. indicate polar parameters of $\A$, while $t_0^A,\phi_1^A$ etc. the corresponding ones of $\Ac$, and similarly for $\B,\Bc,\D$ and $\Dc$). A polar quantity of the basic layer has no superscript.

The polar components of $\A,\B$ and $\D$ are given by
\begin{equation}
 \label{eq:compApolar}
 \A\ \ \rightarrow\ \
 \left\{
 \begin{split}
 &T_0^A=\frac{1}{h}\sum_{k=1}^n{T_0}_k(z_k-z_{k-1}),\\
 &T_1^A=\frac{1}{h}\sum_{k=1}^n{T_1}_k(z_k-z_{k-1}),\\
&R_0^Ae^{4i\Phi_0^A}=\frac{1}{h}\sum_{k=1}^n{R_0}_ke^{4i({\Phi_0}_k+\delta_k)}(z_k-z_{k-1}),\\
&R_1^Ae^{2i\Phi_1^A}=\frac{1}{h}\sum_{k=1}^n{R_1}_ke^{2i({\Phi_1}_k+\delta_k)}(z_k-z_{k-1}).
 \end{split}
 \right.
 \end{equation}
\begin{equation}
 \label{eq:compBpolar}
 \B\ \ \rightarrow\ \
 \left\{
 \begin{split}
 &T_0^B=\frac{1}{h^2}\sum_{k=1}^n{T_0}_k(z_k^2-z_{k-1}^2),\\
 &T_1^B=\frac{1}{h^2}\sum_{k=1}^n{T_1}_k(z_k^2-z_{k-1}^2),\\
&R_0^Be^{4i\Phi_0^B}=\frac{1}{h^2}\sum_{k=1}^n{R_0}_ke^{4i({\Phi_0}_k+\delta_k)}(z_k^2-z_{k-1}^2),\\
&R_1^Be^{2i\Phi_1^B}=\frac{1}{h^2}\sum_{k=1}^n{R_1}_ke^{2i({\Phi_1}_k+\delta_k)}(z_k^2-z_{k-1}^2).
 \end{split}
 \right.
 \end{equation}
\begin{equation}
 \label{eq:compDpolar}
 \D\ \ \rightarrow\ \
 \left\{
 \begin{split}
 &T_0^D=\frac{4}{h^3}\sum_{k=1}^n{T_0}_k(z_k^3-z_{k-1}^3),\\
 &T_1^D=\frac{4}{h^3}\sum_{k=1}^n{T_1}_k(z_k^3-z_{k-1}^3),\\
&R_0^De^{4i\Phi_0^D}=\frac{4}{h^3}\sum_{k=1}^n{R_0}_ke^{4i({\Phi_0}_k+\delta_k)}(z_k^3-z_{k-1}^3),\\
&R_1^De^{2i\Phi_1^D}=\frac{4}{h^3}\sum_{k=1}^n{R_1}_ke^{2i({\Phi_1}_k+\delta_k)}(z_k^3-z_{k-1}^3).
 \end{split}
 \right.
 \end{equation}
It is worth noting that the isotropic and anisotropic parts of all the tensors remain {\it separated} in the homogenization of the polar parameters, for all the tensors. Moreover, {\it special orthotropies are preserved}: 
 \begin{equation}
 \begin{split}
 &{R_0}_k=0\ \forall k\ \Rightarrow\ R_0^A=R_0^B=R_0^D=0,\\
 &{R_1}_k=0\ \forall k\ \Rightarrow\ R_1^A=R_1^B=R_1^D=0.
 \end{split}
 \end{equation}
 In particular, for  a laminate  composed by layers that are all square symmetric, though different,  $\A,\ \B$ and $\D$ are square symmetric too, no matter  the stacking sequence and  the layer orientations. The same is true for $R_0-$orthotropic plies, while ordinary orthotropy is not preserved through the homogenization process. 
 
For what concerns a second-rank symmetric tensor {\bf L}, in the polar formalism it is
\be
\label{eq:mohr}
\begin{split}
&\gr{L}_{1}(\theta)=T+R\cos2(\varPhi-\theta),\\
&\gr{L}_{2}(\theta)=T-R\cos2(\varPhi-\theta),\\
&\gr{L}_{6}(\theta)=\sqrt{2}R\sin2(\varPhi-\theta),
\end{split}
\ee
with $T, R$ two invariants, representing respectively the {\it isotropic} and the {\it anisotropic} phases of {\bf L}; $\Phi$ is an angle determined by the choice of the frame. In the following, we will indicate by 
$T^U,R^U,\Phi^U$ the polar components of $\U$, by $t^u,r^u,\phi^u$ those of $\u$ and similarly for $\V$, $\W$ and $\w$, while $t_1^v,r_1^v,\phi_1^v$ will denote the polar components of $\v_1$ and $t_2^v,r_2^v,\phi_2^v$ those of $\v_2$.

The polar parameters of $\U,\V,\W$ are
\be
\gr{U}\ \rightarrow\ \left\{\begin{split}&T^U=\frac{1}{h}\sum_{k=1}^n{T_\gamma}_k(z_k-z_{k-1}),
\\&R^Ue^{2i\Phi^U}=\frac{1}{h}\sum_{k=1}^n{R_\gamma}_k e^{2i({\Phi_\gamma}_k+\delta_k)}(z_k-z_{k-1}),\end{split}\right.
\ee
\be
\gr{V}\ \rightarrow\ \left\{\begin{split}&T^V=\frac{1}{h^2}\sum_{k=1}^n{T_\gamma}_k(z_k^2-z_{k-1}^2),
\\&R^Ve^{2i\Phi^V}=\frac{1}{h^2}\sum_{k=1}^n{R_\gamma}_ke^{2i({\Phi_\gamma}_k+\delta_k)}(z_k^2-z_{k-1}^2),\end{split}\right.
\ee
\be
\gr{W}\ \rightarrow\ \left\{\begin{split}&T^W=\frac{4}{h^3}\sum_{k=1}^n{T_\gamma}_k(z_k^3-z_{k-1}^3),
\\&R^We^{2i\Phi^W}=\frac{4}{h^3}\sum_{k=1}^n{R_\gamma}_ke^{2i({\Phi_\gamma}_k+\delta_k)}(z_k^3-z_{k-1}^3),\end{split}\right.
\ee 
 where $T_\gamma,R_\gamma,\Phi_\gamma$ are the polar parameters of tensor $\bgamma(\delta_k)$.

For  laminates composed of identical plies, a very important case for applications, the  relations for $\A,\B,\D,$ become
 \begin{equation}
  \label{eq:compABDpolareq}
  \begin{split}
&\mathbb{A}\rightarrow\left\{
\begin{array}{l}
T_0^A=T_0,\\
T_1^A=T_1,\\
R_0^Ae^{4i\varPhi_0^A}={R_0e^{4i\varPhi_0}}{\left(\xi_1+i\xi_2\right)},\\
R_1^Ae^{2i\varPhi_1^A}={R_1e^{2i\varPhi_1}}{\left(\xi_3+i\xi_4\right)};
\end{array}\right.\\
&\mathbb{B}\rightarrow\left\{
\begin{array}{l}
T_0^B=0,\\
T_1^B=0,\\
R_0^Be^{4i\varPhi_0^B}={R_0e^{4i\varPhi_0}}{\left(\xi_5+i\xi_6\right)},\\
R_1^Be^{2i\varPhi_1^B}={R_1e^{2i\varPhi_1}}{\left(\xi_7+i\xi_8\right)};
\end{array}\right.\\
&\mathbb{D}\rightarrow\left\{
\begin{array}{l}
T_0^D=T_0,\\
T_1^D=T_1,\\
R_0^De^{4i\varPhi_0^D}={R_0e^{4i\varPhi_0}}(\xi_{9}+i\xi_{10}),\\
R_1^De^{2i\varPhi_1^D}={R_1e^{2i\varPhi_1}}(\xi_{11}+i\xi_{12}),
\end{array}\right.
\end{split}
\end{equation}
and those for $\U,\V$ and $\W$
\be
\label{eq:UVWpolaridlayers}
\begin{split}
&\gr{U}\ \rightarrow\ \left\{\begin{split}&T^U=T_\gamma,\\&R^Ue^{2i\Phi^U}=R_\gamma e^{2i\Phi_\gamma}(\xi_3+i\xi_4),\end{split}\right.\\
&\gr{V}\ \rightarrow\ \left\{\begin{split}&T^V=0,\\&R^Ve^{2i\Phi^V}=R_\gamma e^{2i\Phi_\gamma}(\xi_7+i\xi_8);\end{split}\right.\\
&\gr{W}\ \rightarrow\ \left\{\begin{split}&T^W=T_\gamma,\\&R^We^{2i\Phi^W}=R_\gamma e^{2i\Phi_\gamma}(\xi_{11}+i\xi_{12}).\end{split}\right.
\end{split}
\ee
In the above equations, the quantities $\xi_i,\ i=1,...,12,$ are the so-called {\it lamination parameters} \cite{TsaiHahn,vannucci_libro},   quantities accounting for the geometry of the stacking sequence (i.e. for orientations and position of the layers on the stack):
\begin{equation}
\label{eq:laminationparameters}
\begin{aligned}
&{{{{\xi}}_{1}{{+}}{i}{{\xi}}_{2}}{{=}}\mathop{\sum}\limits_{{k}{{=}}{1}}\limits^{n}a_k\hspace{0.33em}{{e}^{4i{{\delta}}_{k}}}},\hspace{1cm}
& &{{{{\xi}}_{3}{{+}}{i}{{\xi}}_{4}}{{=}}\mathop{\sum}\limits_{{k}{{=}}{1}}\limits^{n}a_k\hspace{0.33em}{{e}^{2i{{\delta}}_{k}}}},\\
&{{{{\xi}}_{5}{{+}}{i}{{\xi}}_{6}}{{=}}\mathop{\sum}\limits_{{k}{{=}}{1}}\limits^{n}{{b}_{k}\hspace{0.33em}{e}^{4i{{\delta}}_{k}}}},\hspace{1cm}
& &{{{{\xi}}_{7}{{+}}{i}{{\xi}}_{8}}{{=}}\mathop{\sum}\limits_{{k}{{=}}{1}}\limits^{n}{{b}_{k}\hspace{0.33em}{e}^{2i{{\delta}}_{k}}}},\\
&{{{{\xi}}_{9}{{+}}{i}{{\xi}}_{10}}{{=}}\mathop{\sum}\limits_{{k}{{=}}{1}}\limits^{n}{{d}_{k}\hspace{0.33em}{e}^{4i{{\delta}}_{k}}}},\hspace{1cm}
& &{{{{\xi}}_{11}{{+}}{i}{{\xi}}_{12}}{{=}}\mathop{\sum}\limits_{{k}{{=}}{1}}\limits^{n}{{d}_{k}\hspace{0.33em}{e}^{2i{{\delta}}_{k}}}},\\%
\end{aligned}
\end{equation}
with
\begin{equation}
\label{eq:coefABCD}
a_k=\frac{1}{n},\ b_k=\frac{1}{n^2}(2k-n-1),\
d_k=\frac{1}{n^3}\left[12k(k-n-1)+4+3n(n+2)\right].
\end{equation}

So, through the polar formalism we can see that:
\begin{itemize}
\item the isotropic part of $\A$ and $\D$ is equal to that of the basic layer: $T_0^A=T_0^D=T_0,T_1^A=T_1^D=T_1$;
\item in the same way, the isotropic part of $\U$ and $\W$ is equal to that of $\bgamma:T^U=T^W=T_\gamma$;
\item $\B$ is exclusively anisotropic, i.e. its isotropic part is null, $T_0^B=T_1^B=0$, so its average on $\{0,2\pi\}$  is zero;
\item $\V$ is exclusively anisotropic, i.e. its isotropic part is null, $T^V=0$, so also its average on  $\{0,2\pi\}$  is zero.
\end{itemize}

Tensors $\B$ and $\Bc$ are indeed special operators. On the one hand, for laminates made of identical plies  $T_0^B=T_1^B$, which gives that $\B$ is a {\it rari-constant tensor}, \cite{vannucci2016}: it has also the Cauchy-Poisson symmetries, \cite{pearson,love,benvenuto91}, besides the minor and major ones. In practice,  $\B$ is completely symmetric with respect to any permutation of the indices, i.e. it is also
\be
\B_{1122}=\B_{1212}\ \rightarrow\ \B_{12}=\frac{\B_{66}}{2}.
\ee
However, for hybrid laminates in general $T_0^B\neq T_1^B$, so $\B$ is not rari-constant for hybrid coupled laminates. 

Concerning $\Bc$, we have already noticed that  $\Bc\neq \Bc^\top$, i.e. it does not possess the major symmetries. Contrarily to $\B$, for identical layers laminates, that has a greater number of symmetries with respect to a usual elastic tensor, $\Bc$ has a lacking of  symmetries. Such a special  elastic tensor has been studied in the framework of the polar formalism, \cite{vannucci10ijss}. In particular, it has been shown that a tensor of this type has nine independent  components:
\begin{equation}
\label{eq:bspecial}
\begin{split}
&\Bc_{11}(\theta)=t_0^B\hspace{-1mm}+\hspace{-1mm}2t_1^B\hspace{-1mm}+\hspace{-1mm}r_0^B\cos4(\phi_0^B\hspace{-1mm}-\hspace{-1mm}\theta)+2r_1^B\cos2(\phi_1^B\hspace{-1mm}-\hspace{-1mm}\theta)+\\
&\hspace{10mm}+2r_2^B\cos2(\phi_2^B-\theta),\\
&\Bc_{16}(\theta)=\sqrt{2}[-t_3^B+\hspace{-1mm}r_0^B\sin4(\phi_0^B\hspace{-1mm}-\hspace{-1mm}\theta)+\hspace{-1mm}2r_2^B\sin2(\phi_2^B-\hspace{-1mm}\theta)],\\
&\Bc_{12}(\theta)=-t_0^B\hspace{-1mm}+\hspace{-1mm}2t_1^B\hspace{-1mm}-\hspace{-1mm}r_0^B\cos4(\phi_0^B\hspace{-1mm}-\hspace{-1mm}\theta)\hspace{-1mm}+\hspace{-1mm}2r_1^B\cos2(\phi_1^B\hspace{-1mm}-\hspace{-1mm}\theta)-\\
&\hspace{10mm}-2r_2^B\cos2(\phi_2^B-\theta),\\
&\Bc_{61}(\theta)=\sqrt{2}[t_3^B+\hspace{-1mm}r_0^B\sin4(\phi_0^B-\hspace{-1mm}\theta)+2r_1^B\sin2(\phi_1^B-\hspace{-1mm}\theta)],\\
&\Bc_{66}(\theta)=2[t_0^B-r_0^B\cos4(\phi_0^B-\theta)],\\
&\Bc_{62}(\theta)=\sqrt{2}[-t_3^B\hspace{-1mm}-\hspace{-1mm}r_0^B\sin4(\phi_0^B\hspace{-1mm}-\hspace{-1mm}\theta)+\hspace{-1mm}2r_1^B\sin2(\phi_1^B-\theta)],\\
&\Bc_{21}(\theta)=-t_0^B\hspace{-1mm}+\hspace{-1mm}2t_1^B\hspace{-1mm}-\hspace{-1mm}r_0^B\cos4(\phi_0^B\hspace{-1mm}-\hspace{-1mm}\theta)\hspace{-1mm}-\hspace{-1mm}2r_1^B\cos2(\phi_1^B\hspace{-1mm}-\hspace{-1mm}\theta)+\\
&\hspace{10mm}+2r_2^B\cos2(\phi_2^B-\theta),\\
&\Bc_{26}(\theta)=\sqrt{2}[t_3^B-r_0^B\sin4(\phi_0^B-\hspace{-1mm}\theta)+2r_2^B\sin2(\phi_2^B-\hspace{-1mm}\theta)],\\
&\Bc_{22}(\theta)=t_0^B+2t_1^B+r_0^B\cos4(\phi_0^B\hspace{-1mm}-\hspace{-1mm}\theta)\hspace{-1mm}-\hspace{-1mm}2r_1^B\cos2(\phi_1^B\hspace{-1mm}-\hspace{-1mm}\theta)-\\
&\hspace{10mm}-2r_2^B\cos2(\phi_2^B-\theta).\\
\end{split}
\end{equation}
In practice, for this tensor, three other polar parameters appear: $t_3^B,r_2^B$ and $\phi_2^B$, as well as three more independent Cartesian components: $\Bc_{21},\Bc_{61}$ and $\Bc_{62}$.

\section{Recall of known results}
Some general results, found through the polar method, are recalled hereafter, cf. \cite{vannucci_libro}. We just recall that a laminate is said to be {\it thermally uncoupled} when
\be
\V=\v_1=\v_2=\mathbf{O}.
\ee
\subsection{Hybrid laminates}
For hybrid laminates elastic uncoupling is {\it not sufficient} to ensure also thermoelastic uncoupling: $\B=\mathbb{O}\nRightarrow \V=\mathbf{O}$: they can exist laminates that are elastically uncoupled but thermally coupled and in this case, putting $\B=\mathbb{O}$ into eqs. (\ref{eq:inversetensorsABD}) and (\ref{eq:tenstherminverses}), we get
\be
\u=\A^{-1}\U,\ \v_1=\frac{6}{h}\D^{-1}\V,\ \v_2=\frac{h}{2}\A^{-1}\V,\ \w=\D^{-1}\W.
\ee
This shows that a hybrid laminate is thermally uncoupled if and only if 
\be
\B=\mathbb{O}\ and\ \V=\mathbf{O}.
\ee
No other general conclusions can be found for hybrid laminates.

\subsection{Laminates made of identical layers}
The situation is different for laminates made of identical layers: from eqs. (\ref{eq:compABDpolareq}),  (\ref{eq:UVWpolaridlayers}),  (\ref{eq:inversetensorsABD}) and (\ref{eq:tenstherminverses}) we have that
\be
\B=\mathbb{O}\Rightarrow\xi_7+i\xi_8=0\Rightarrow\V=\v_1=\v_2=\mathbf{O}:
\ee
for a laminate composed of identical layers, $\B=\mathbb{O}$ is a sufficient condition for thermal uncoupling; we remark that for these laminates $R_1^B=0\iff R^V=0$. 

The  converse is not true: $\V=\mathbf{O}\nRightarrow\B=\mathbb{O}$, because $R_0^B$ can be different from zero. In this particular case,
\be
\begin{split}
&\gr{u}=\mathcal{A}\gr{U}=(\A-3\B\D^{-1}\B)^{-1}\gr{U},\\ 
&\gr{v}_1=\frac{2}{h}\mathcal{B}^\top\gr{U}=-\frac{6}{h}(\D-3\B\A^{-1}\B)^{-1}\B\A^{-1}\gr{U},\\
&\gr{v}_2=\frac{h}{6}\mathcal{B}\gr{W}=-\frac{h}{2}(\A-3\B\D^{-1}\B)^{-1}\B\D^{-1}\gr{W},\\ 
&\gr{w}=\mathcal{D}\gr{W}=(\D-3\B\A^{-1}\B)^{-1}\gr{W}.
\end{split}
\ee
i.e. $\gr{v}_1$ and $\gr{v}_2$ are not necessarily null: a change of temperature $t$ makes the plate warp but does not produce bending moments and a gradient $\nabla t$ stretches the plate but does not give rise to membrane forces. In small words, it is in principle possible to obtain laminates that, for what regards coupled effects, are sensitive to temperature changes only for what concerns the deformation, but insensitive for internal actions. For such laminates $\B$ is square symmetric ($R_1^B=0$).

\subsection{Quasi homogeneous coupled laminates}
\label{sec:TQHCL}
The concept of {\it quasi homogeneity} for laminates has been introduced  by G. Verchery in 1988, \cite{Kandil}: defined the {\it homogeneity tensor} 
 \begin{equation}
 \label{eq:quasihomogeneitydefinition}
 \C:=\A-\D,
 \end{equation}
a laminate is said to be {\it coupled quasi-homogeneous} $\iff$
 \begin{equation}
 \label{eq:quasihomogene}
\C=\mathbb{O}.
 \end{equation}
For these laminates, the elastic behavior is the same in bending and in extension. By eqs. (\ref{eq:compABD})  we get 
  \begin{equation}
   \label{eq:compC}
 \C=\frac{1}{h^3}\sum_{k=1}^n\left[h^2(z_k-z_{k-1})-4({z_k^3}-{z_{k-1}^3})\right]\Q_k(\delta_k),
\end{equation}
that for the case of identical layers becomes
  \begin{equation}
  \label{eq:abcdeqply}
\C=\sum_{k=1}^nc_k\Q(\delta_k),
 \end{equation}
where
\begin{equation}
\label{eq:coeffC}
c_k=a_k-d_k=\frac{1}{n^3}[-2(6k^2+n^2)+6n(2k-1)+4(3k-1)].
\end{equation}
By the polar method, in this last case we get also
\begin{equation}
  \label{eq:compCpolareq}
 \C \rightarrow
 \left\{
 \begin{split}
 &T_0^C=0,\\
 &T_1^C=0,\\
&R_0^Ce^{4i\Phi_0^C}={R_0}e^{4i{\Phi_0}}{\left(\xi_9+i\xi_{10}\right)},\\
&R_1^Ce^{2i\Phi_1^C}={R_1}e^{2i{\Phi_1}}{\left(\xi_{11}+i\xi_{12}\right)},
 \end{split}
 \right.
 \end{equation}
with the lamination parameters
\begin{equation}
\label{eq:laminationparC}
\begin{aligned}
&{{{{\xi}}_{13}{{+}}{i}{{\xi}}_{14}}{{=}}\mathop{\sum}\limits_{{j}{{=}}{1}}\limits^{n}{{c}_{j}\hspace{0.33em}{e}^{4i{{\delta}}_{j}}}},\hspace{1cm}
& &{{{{\xi}}_{15}{{+}}{i}{{\xi}}_{16}}{{=}}\mathop{\sum}\limits_{{j}{{=}}{1}}\limits^{n}{{c}_{j}\hspace{0.33em}{e}^{2i{{\delta}}_{j}}}}.
\end{aligned}
\end{equation}

The concept of coupled quasi homogeneity has been extended to thermal properties in \cite{vannucci_libro}: a laminate is  {\it thermally quasi-homogeneous coupled} if, in addition to $\C=\mathbb{O}$, it is also
\be
\gr{Y}=\gr{O},
\ee
where 
\be
\label{eq:quasihomogeneitydefinitionthermal}
\gr{Y}:=\gr{U}-\gr{W}
\ee
is the {\it thermal homogeneity tensor}. From eq. (\ref{eq:tensorthemoelast}) we get
\be
\gr{Y}=\frac{1}{h^3}\sum_{k=1}^n[h^2(z_k-z_{k-1})-4(z_k^3-z_{k-1}^3)]\bgamma_k(\delta_k),
\ee
and for the case of identical layers
\be
\gr{Y}=\sum_{k=1}^nc_k\bgamma_k(\delta_k),
\ee
with the same coefficients $c_k$  given by eq. (\ref{eq:coeffC}). In this case,  the polar components of $\gr{Y}$ are
\be
\gr{Y}\ \rightarrow\ \left\{\begin{split}&T^Y=0,\\&R^Ye^{2i\Phi^Y}=R_\gamma e^{2i\Phi_\gamma}(\xi_{15}+i\xi_{16}).\end{split}\right.
\ee
Because for laminates made of identical plies
\be
\C=\mathbb{O} \Rightarrow \xi_{15}+i\xi_{16}=0\Rightarrow \gr{Y}=\gr{O},
\ee
we see that {\it coupled quasi-homogeneity implies  thermal coupled quasi-homogeneity}. However, the converse is not true, because $R_0^C$ can be different from zero: $\gr{Y}=\gr{O}\nRightarrow\C=\mathbb{O}$. This means that we can find laminates with an equal thermo-elastic anisotropic behavior in extension and in bending, i.e. with the coefficients of thermal expansion and curvature equal together for all the directions, but with elastic behaviors in membrane and bending that are different.

In this paper, we are interested  in  {\it  thermally quasi homogeneous coupled laminates} (TQHCL), i.e. having $\gr{Y}=\gr{O}$ but at least one among $\v_1$ and $\v_2$ different from zero (the condition $\V\neq\gr{O}$ is not sufficient, because laminates having $\V=\gr{O}$ but  $\v_1\neq\gr{O}$ or $\v_2\neq\gr{O}$ can exist, see below). 

We can prove directly an interesting result: we already know that $\A=\D\Rightarrow\Ac=\Dc\ and\ \Bc=\Bc^\top$, see \cite{vannucci23a}; moreover, see above, $\C=\A-\D=\mathbb{O}\Rightarrow\gr{Y}=\U-\W=\gr{O}$, so, from eq. (\ref{eq:tenstherminverses}) we get
\be
\u=\Ac\U+\Bc\V=\Dc\W+\Bc^\top\V=\w,
\ee
i.e. the condition $\A=\D$ implies not only $\Ac=\Dc$, but also $\U=\W$ and $\u=\w$: quasi homogeneity is a property that holds for stiffness and compliance, for both the elastic and thermal responses, and for coupled or uncoupled laminates.

\section{Relations between the elastic and thermoelastic  tensors}
The polar formalism allows for making explicit some relations among the elastic and the thermoelastic tensors, of a laminate, both for the stiffness and the compliance tensors. These relations are introduced and examined in the following Sections. From now on, we refer exclusively to laminates made of identical plies, all of them orthotropic and with a material frame fixed by the choice $\Phi_1=0$:
\begin{equation}
\label{eq:Q}
\begin{split}
&\mathbb{Q}_{11}{(}\theta{)}{=}{T}_{0}{+}{2}{T}_{1}{+}(-1)^k{R}_{0}\cos{4}\theta{+}{4}{R}_{1}\cos{2}\theta,\\
&\mathbb{Q}_{12}{(}\theta{)}{=}{-}{T}_{0}{+}{2}{T}_{1}{-}(-1)^k{R}_{0}\cos{4}\theta,\\
&\mathbb{Q}_{16}(\theta)=-\sqrt{2}\left[(-1)^kR_0\sin4\theta+2R_1\sin2\theta\right],\\
&\mathbb{Q}_{22}{(}\theta{)}{=}{T}_{0}{+}{2}{T}_{1}{+}(-1)^k{R}_{0}\cos{4}\theta{-}{4}{R}_{1}\cos{2}\theta,\\
&\mathbb{Q}_{26}{(}\theta{)}{=}\sqrt{2}\left[(-1)^k{R}_{0}\sin{4}\theta{-}{2}{R}_{1}\sin{2}\theta\right],\\
&\mathbb{Q}_{66}{(}\theta{)}{=}2\left[{T}_{0}{-}(-1)^k{R}_{0}\cos{4}\theta\right].
\end{split}
\end{equation}
Moreover,  the tensor of thermal expansion coefficients $\balpha$ is also orthotropic (or isotropic) with the same orthotropy axes of $\Q$, so that
\be
\Phi_\gamma=\Phi_1+\lambda\frac{\pi}{2},\ \ \lambda\in\{0,1\}.
\ee
This hypothesis actually corresponds to the physical reality of composite layers.

\subsection{Relations concerning the stiffness tensors}
We first consider the relations existing among the elastic and thermoelastic stiffness tensors.

\subsubsection{A physical constant}
From eq. (\ref{eq:compABDpolareq})  we get the following relations
\be
R_1^A=R_1\sqrt{\xi_3^2+\xi_4^2},\ \ R_1^B=R_1\sqrt{\xi_7^2+\xi_8^2},\ \ R_1^D=R_1\sqrt{\xi_{11}^2+\xi_{12}^2},
\ee
while by eq.  (\ref{eq:UVWpolaridlayers}) 
\be
R^U=R_\gamma\sqrt{\xi_3^2+\xi_4^2},\ \ R^V=R_\gamma\sqrt{\xi_7^2+\xi_8^2},\ \ R^W=R_\gamma\sqrt{\xi_{11}^2+\xi_{12}^2}.
\ee
Making the ratio of the respective terms, we have then
\be
\frac{R^U}{R_1^A}=\frac{R^V}{R_1^B}=\frac{R^W}{R_1^D}=\frac{R_\gamma}{R_1}:=\rho;
\ee
it exists hence a physical constant, $\rho$, which is a ratio between the corresponding couples of elastic and thermoelastic tensors for the layer and for the laminate, i.e. $\Q$ and $\bgamma$, $\A$ and $\U$, $\B$ and $\V$, $\D$ and $\W$. The choice of a material fixes the value of $\rho$, so the same ratio also for the stiffness tensors of the laminate.

\subsubsection{Tensor orientations}
\begin{figure}[h]
\center
\includegraphics[width=.5\textwidth]{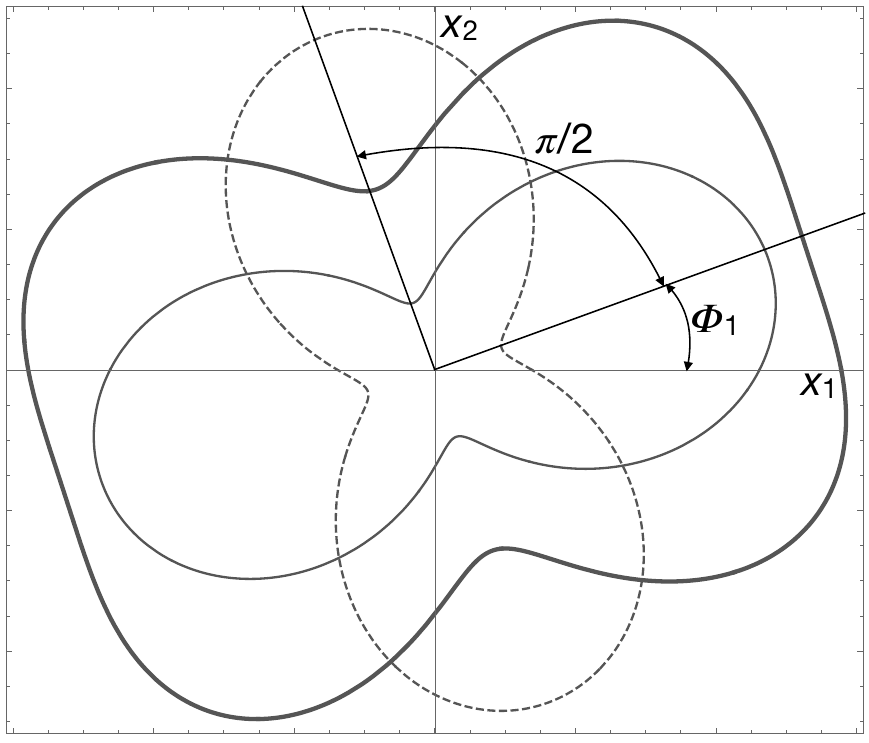}
\caption{Possible situations of the components $\A_{11}$, thick line, and $U_1$, thin and dashed lines, for a laminate.}
\label{fig:2}
\end{figure}

Still from eq. (\ref{eq:compABDpolareq}) we can calculate the polar angles of $\A,\B$ and $\D$:
\be
\begin{array}{ll}
\Phi_0^A=\dfrac{1}{4}\arctan\dfrac{\xi_2}{\xi_1}+k\dfrac{\pi}{4},& \Phi_1^A=\dfrac{1}{2}\arctan\dfrac{\xi_4}{\xi_3},\bigskip\\
\Phi_0^B=\dfrac{1}{4}\arctan\dfrac{\xi_6}{\xi_5}+k\dfrac{\pi}{4},& \Phi_1^B=\dfrac{1}{2}\arctan\dfrac{\xi_8}{\xi_7},\bigskip\\
\Phi_0^D=\dfrac{1}{4}\arctan\dfrac{\xi_{10}}{\xi_9}+k\dfrac{\pi}{4},& \Phi_1^D=\dfrac{1}{2}\arctan\dfrac{\xi_{12}}{\xi_{11}},
\end{array}
\ee 
and by eq. (\ref{eq:UVWpolaridlayers}) those of $\U,\V$ and $\W$:
\be
\Phi^U=\frac{1}{2}\arctan\frac{\xi_4}{\xi_3}+\lambda\frac{\pi}{2},\ 
\Phi^V=\frac{1}{2}\arctan\frac{\xi_8}{\xi_7}+\lambda\frac{\pi}{2},\ 
\Phi^W=\frac{1}{2}\arctan\frac{\xi_{12}}{\xi_{11}}+\lambda\frac{\pi}{2}.
\ee
Comparing these equations, we get immediately
\be
\label{eq:angoliUVW}
\Phi^U=\Phi_1^A+\lambda\frac{\pi}{2},\  \Phi^V=\Phi_1^B+\lambda\frac{\pi}{2},\ \Phi^W=\Phi_1^D+\lambda\frac{\pi}{2}:
\ee
each stiffness thermal tensor, $\U,\V$ or $\W$, is aligned or rotated of $\pi/2$ with respect to the corresponding elastic stiffness tensor, cf. Fig. \ref{fig:2}. To remark that $\U,\V$ and $\W$, as planar 2nd-rank symmetric tensors, are either orthotropic or isotropic, while $\A,\B$ and $\D$ can be completely anisotropic: the above relations hold in all the cases (if a thermal tensor is isotropic, it can be conventionally assumed that its polar angle is equal to that of the corresponding elastic tensor). 

\subsubsection{Relations among symmetries}
The existence of the ratio $\rho$ gives that
\begin{enumerate}[i.]
\item the anisotropic phase of the thermoelastic tensors is fixed by the $R_1$-phase of the elastic ones:
\be
R^U=\rho\ R_1^A,\ \ R^V=\rho\ R_1^B,\ \ R^W=\rho\ R_1^D;
\ee
\item the square symmetry of $\A,\D$ implies the isotropy of $\U,\W$ respectively:
\be
R_1^A=0\Rightarrow R^U=0,\ \ R_1^D=0\Rightarrow R^W=0;
\ee
\item {\it a fortiori}, this is true also for the case of $\A$ and $\D$ isotropic;
\item the square symmetry of $\B$ implies the vanishing of $\V$:
\be
R_1^B=0\Rightarrow R^V=0\Rightarrow \V=\gr{O};
\ee
\item laminates composed of square symmetric layers have $R_\gamma=R_1=0\Rightarrow\rho$ is not defined; however, from eqs. (\ref{eq:compABDpolareq}) and (\ref{eq:UVWpolaridlayers}), we get again that $R_1^A=R_1^B=R_1^D=R^U=R^V=R^W=0$, no matter of the orientation angle of the plies, for both the cases of coupled, $\B\neq\mathbb{O}$, or uncoupled, $\B=\mathbb{O}$, laminates.
\end{enumerate}

\subsection{Relations concerning the compliance tensors}
\begin{figure}[h]
\center
\includegraphics[width=.6\textwidth]{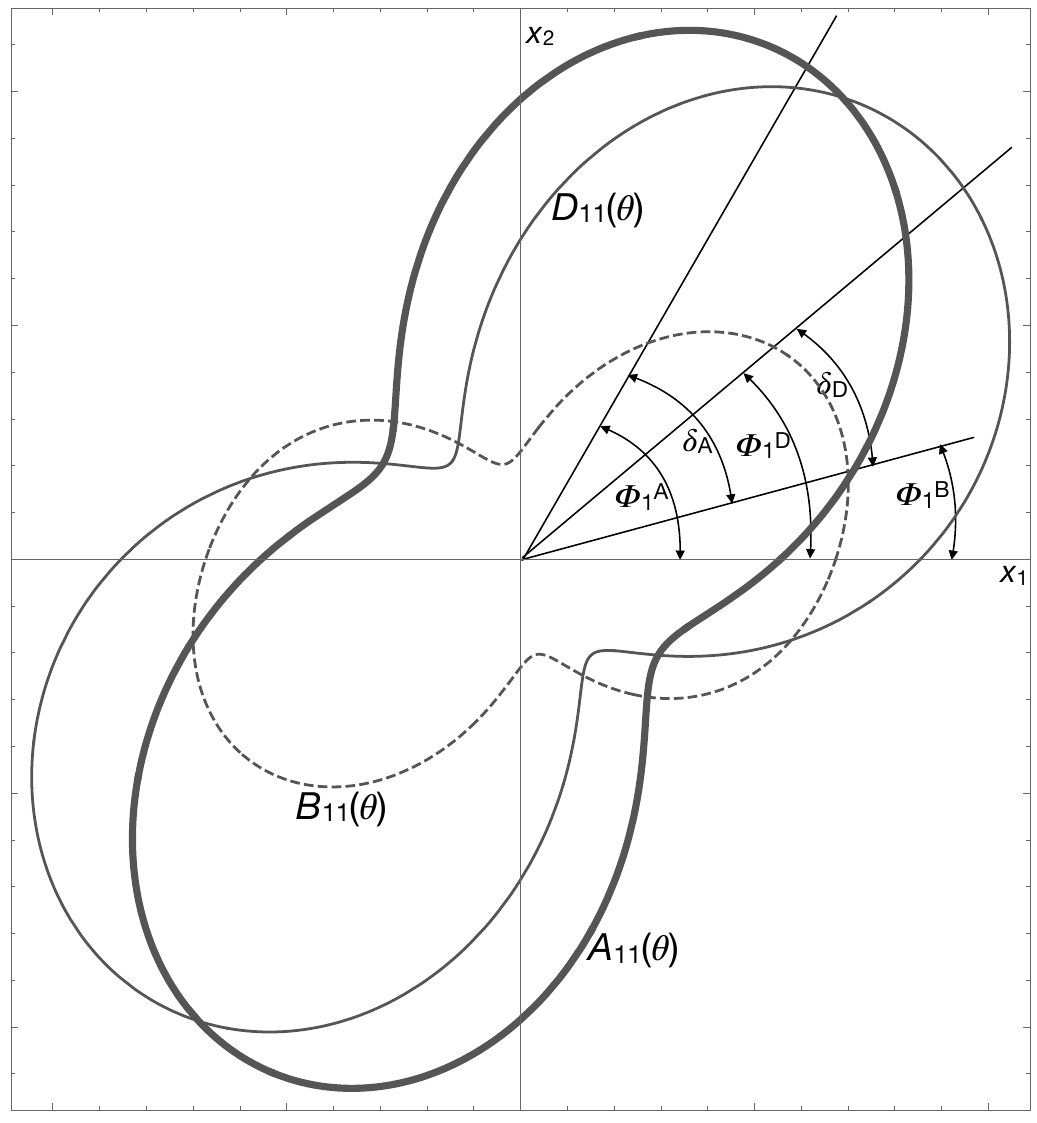}
\caption{The polar and shift angles shown for a generic material.}
\label{fig:3}
\end{figure}

We consider now the relations between the stiffness tensors, $\A,\B,\D,\U,\V,\W$, and the compliance thermoelastic tensors $\u,\v_1,\v_2,\w$. The relations among $\A,\B,\D$ and $\Ac,\Bc,\Dc$ for the case of coupled laminates have been already studied in \cite{vannucci23a}, the results are used below. In the following, we will use the angle differences
\be
\Phi_A=\Phi_0^A-\Phi_1^A,\ \ \Phi_B=\Phi_0^B-\Phi_1^B,\ \ \Phi_D=\Phi_0^D-\Phi_1^D,
\ee
that are tensor invariants, and the {\it shift angles}
\be
\delta_A=\Phi_1^A-\Phi_1^B,\ \ \delta_D=\Phi_1^D-\Phi_1^B,
\ee
that define the relative position of tensors $\A$ and $\D$ with respect to tensor $\B$, see Fig. \ref{fig:3}. Then, all the equations are written putting $\theta=\Phi_1^B$ in eqs. (\ref{eq:mohr4}) and (\ref{eq:mohr}) for all the tensors $\A,\B,\D,\U,\V$ and $\W$ or, which is the same, we determine the $x_1$ axis fixing $\Phi_1^B=0$; finally, we use eq. (\ref{eq:angoliUVW}), so that, with the last assumption,
\be
\Phi^U=\delta_A+\lambda\frac{\pi}{2},\ \ \Phi^V=\lambda\frac{\pi}{2},\ \ \Phi^W=\delta_D+\lambda\frac{\pi}{2}.
\ee

The general case of completely anisotropic and different tensors $\A,\B$ and $\D$ has been considered but it is not reported below, because the results are very cumbersome and in practice too long  to be rewritten; moreover, such a general case is not interesting for practical applications, when some type of elastic symmetry, namely orthotropy, is always prescribed. So, we present hereafter some cases characterized by the existence of some elastic symmetries for $\A,\B$ or/and $\D$ and  consider the form of the corresponding compliance thermoelastic tensors $\u,\v_1,\v_2,\w$.

Before going on, we recall that with the above assumptions and through eqs. (\ref{eq:mohr}) and (\ref{eq:UVWpolaridlayers}), it is
\be
\label{eq:formuleUVW}
\begin{split}
&\U=\left\{\begin{array}{c}T_\gamma+(-1)^\lambda\rho\ R_1^A\cos2\delta_A\\T_\gamma-(-1)^\lambda\rho\ R_1^A\cos2\delta_A\\(-1)^\lambda\rho\ R_1^A\sin2\delta_A\end{array}\right\},\\
&\V=\left\{\begin{array}{c}(-1)^\lambda\rho\ R_1^B\\-(-1)^\lambda\rho\ R_1^B\\0\end{array}\right\},\\
&\W=\left\{\begin{array}{c}T_\gamma+(-1)^\lambda\rho\ R_1^D\cos2\delta_D\\T_\gamma-(-1)^\lambda\rho\ R_1^D\cos2\delta_D\\(-1)^\lambda\rho\ R_1^D\sin2\delta_D\end{array}\right\},
\end{split}
\ee
and
\be
\begin{split}
&\u=\left\{\begin{array}{c}t^u+r^u\cos2\phi^u\\t^u-r^u\cos2\phi^u\\r^u\sin2\phi^u\end{array}\right\},\ \
\v_1=\left\{\begin{array}{c}t_1^v+r_1^v\cos2\phi_1^v\\t_1^v-r_1^v\cos2\phi_1^v\\r_1^v\sin2\phi_1^v\end{array}\right\},\\
&\v_2=\left\{\begin{array}{c}t_2^v+r_2^v\cos2\phi_2^v\\t_2^v-r_2^v\cos2\phi_2^v\\r_2^v\sin2\phi_2^v\end{array}\right\},\ \
\w=\left\{\begin{array}{c}t^w+r^w\cos2\phi^w\\t^w-r^w\cos2\phi^w\\r^w\sin2\phi^w\end{array}\right\}.
\end{split}
\ee
Below, the polar components of the tensors will be specified; through the above formulae, the Cartesian ones can be found immediately. The computation process is as follows: for a particular choice of laminate, characterized by special values of the polar invariants and shift angles determining $\A,\B$ and $\D$, the compliance tensors $\Ac,\Bc$ and $\Dc$ are calculated by eq. (\ref{eq:inversetensorsABD}), while thermal tensors $\U,\V,\W$ by eqs. (\ref{eq:formuleUVW}); then, through eq. (\ref{eq:tenstherminverses}) the Cartesian components of $\u,\v_1,\v_2$ and $\w$ are computed and finally their polar parameters are calculated through te converse of eqs. (\ref{eq:mohr}), \cite{vannucci_libro}. In this way, these last parameters are determined as functions of the stiffness polar parameters of the laminate. The calculations are carried out using a code for analytical computations. The use of the lamination parameters in conjunction with eqs. (\ref{eq:compABDpolareq}) and (\ref{eq:UVWpolaridlayers}) should allow to express the polar parameters of $\u,\v_1,\v_2$ and $\w$ as functions of those of the basic layer and, of course of the $\xi_j, j=1,...,12$; this has also been done, but the results are expressed by so huge formulae that it is impossible to handle and write them and, more important, to extract general considerations from such equations.

\subsubsection{Coupled laminates with $\V=\gr{O}$}
We first consider some cases of laminates having $R_1^B=R^V=0\Rightarrow\B\neq\mathbb{O},\V=\gr{O}$. $\B$ is hence square symmetric and, as specified before, for this class of laminates a change of temperature does not produce internal actions as an effect of coupling, but only deformations.

\begin{enumerate}[i.]
\item $\A\neq\D$ but $k=0$ orthotropic and aligned with $\B\ (\delta_A=\delta_D=0)$: in this case we obtain
\be
\hspace{-3mm}
\begin{array}{c}
{t^u\hspace{-1mm}=\hspace{-.5mm} }\dfrac{\left[\hspace{-.7mm}({R_0^A}\hspace{-1mm}+\hspace{-1mm}{T_0}) \hspace{-1mm}\left(({R_0^D}\hspace{-1mm}+\hspace{-1mm}{T_0}) {T_1}\hspace{-1mm}-\hspace{-1mm}2 {R_1^D}^2\right)\hspace{-1mm}-\hspace{-1mm}3 {R_0^B}^2 {T_1}\right]
{T_\gamma }\hspace{-1mm}+\hspace{-1mm}2 (\hspace{-.7mm}-\hspace{-.7mm}1)^{\lambda } {R_1^A}^2\hspace{-1mm} \left[2 {R_1^D}^2\hspace{-1mm}-\hspace{-1mm}({R_0^D}\hspace{-1mm}+\hspace{-1mm}{T_0}) {T_1}\hspace{-1mm}\right]\hspace{-1mm} \rho }{4 \tau_1},\medskip\\
{r^u= } \left|{\dfrac{{R_1^A} \left[({R_0^D}+{T_0}) {T_1}-2 {R_1^D}^2\right] \left({T_\gamma }-(-1)^{\lambda
} {T_1} \rho \right)}{2\tau_1}}\right|,\medskip\\
\phi^u=0,
\end{array}
\ee
\be
\begin{array}{c}
{t^v_1= }-\dfrac{3}{h}{R_0^B} {R_1^A} {R_1^D}\dfrac{ \left({T_\gamma }-(-1)^{\lambda } {T_1} \rho \right)}{\tau_1},\medskip\\
{r^v_1= }\dfrac{3}{h} {R_0^B} {R_1^A} {T_1}\left|{\dfrac{ {T_\gamma }-(-1)^{\lambda } {T_1} \rho }{\tau_1}}\right|,\medskip\\
\phi^v_1=0,
\end{array}
\ee
\be
\begin{array}{c}
{t^v_2= }-\dfrac{h}{4}{R_0^B} {R_1^A} {R_1^D}\dfrac{  \left({T_\gamma }-(-1)^{\lambda } {T_1} \rho \right)}{\tau_1},\medskip\\
{r^v_2= }\dfrac{h}{4} {R_0^B} {R_1^D} {T_1}\left|{\dfrac{  {T_\gamma }-(-1)^{\lambda } {T_1} \rho
}{\tau_1}}\right|,\medskip\\
\phi^v_1=0,
\end{array}
\ee
\be
\hspace{-3mm}
\begin{array}{c}
{t^w\hspace{-1mm}=\hspace{-.5mm} }\dfrac{\left[\hspace{-.7mm}({R_0^D}\hspace{-1mm}+\hspace{-1mm}{T_0}) \hspace{-1mm}\left(({R_0^A}\hspace{-1mm}+\hspace{-1mm}{T_0}) {T_1}\hspace{-1mm}-\hspace{-1mm}2 {R_1^A}^2\right)\hspace{-1mm}-\hspace{-1mm}3 {R_0^B}^2 {T_1}\right]
{T_\gamma }\hspace{-1mm}+\hspace{-1mm}2 (\hspace{-.7mm}-\hspace{-.7mm}1)^{\lambda } {R_1^A}^2\hspace{-1mm} \left[2 {R_1^A}^2\hspace{-1mm}-\hspace{-1mm}({R_0^A}\hspace{-1mm}+\hspace{-1mm}{T_0}) {T_1}\hspace{-1mm}\right]\hspace{-1mm} \rho }{4 \tau_1},\medskip\\
{r^w= } \left|{\dfrac{{R_1^D} \left[({R_0^A}+{T_0}) {T_1}-2 {R_1^A}^2\right] \left({T_\gamma }-(-1)^{\lambda
} {T_1} \rho \right)}{2\tau_1}}\right|,\medskip\\
\phi^u=0,
\end{array}
\ee
with
\be
\hspace{-5mm}\tau_1\hspace{-.7mm}=2 {R_1^A}^2\hspace{-1mm} \left[2 {R_1^D}^2-({R_0^D}+{T_0}) {T_1}\right]+{T_1} \hspace{-1mm}\left[({R_0^A}+{T_0})\left(({R_0^D}+{T_0}) {T_1}-2 {R_1^D}^2\right)-3 {R_0^B}^2 {T_1}\right]\hspace{-.7mm}.
\ee
Some remarks can be done about these results:
\begin{itemize}
\item physically, it cannot be excluded that the laminate has an {\it anomalous thermal response}, i.e. that, globally, it shrinks for an increase of temperature $t$ or, for what concerns bending, that it takes a curvature on the opposite sign of $\nabla t$; this means that the two isotropic invariants $t^u$ and $t^w$ are not necessarily positive definite, nor the invariant $\tau_1$; {\it a fortiori}, this is true also for $t^v_1$ and $t_2^v$, as  coupling tensors are not definite, cf. \cite{vannucci23a};
\item for tensor $\v_1$, we notice that the ratio between the anisotropic and the isotropic phases is
\be
\left|\frac{r^v_1}{t^v_1}\right|=\frac{T_1}{R_1^D},
\ee
hence it depends exclusively on the material and on the bending stiffness;
\item for tensor $\v_2$, it is
\be
\left|\frac{r^v_2}{t^v_2}\right|=\frac{T_1}{R_1^A},
\ee
so this ratio depends on the material and on the extension stiffness;
\item the ratio between the two isotropic phases of $\v_1$ and $\v_2$ is proportional to the square of the plate's thickness $h$:
\be
\label{eq:ratiotv}
\frac{t^v_2}{t^v_1}=\frac{h^2}{12},
\ee
in particular, it does not depend on the material and laminate stiffnesses;
\item the ratio between the two anisotropic phases of $\v_1$ and $\v_2$ is still proportional to the square of the plate's thickness but it depends also on the membrane and bending stiffnesses:
\be
\label{eq:ratiorv}
\frac{r^v_2}{r^v_1}=\frac{h^2}{12}\frac{R_1^D}{R_1^A};
\ee
\item unlike the stiffness thermal coupling tensor $\V$, that, because the plies are identical, has a null isotropic part, so it has a null mean, the compliance coupling tensors $\v_1$ and $\v_2$ have an isotropic phase; this confirms what already found also for tensor $\Bc$, \cite{vannucci23a};
\item if, in place of $k=0$ it is $k=1$, the results are globally the same, but in the previous formulae some terms change of sign;  we omit to explicitly treat, here and in the following,  the case $k=1$, as it is less important and because it does not enrich substantially the set of possible cases.
\end{itemize}

\item TQHCL orthotropic laminates, i.e. $\A=\D$ and again $\delta_A=\delta_D=0$; in such a case, cf. Section \ref{sec:TQHCL}, we already know that $\u=\w$; basically, the results are like those of the previous case, it is sufficient to put $R_D=R_0^A$ and $R_1^D=R_1^A$; the remarks are the same, unless for the ratio of the anisotropic parts of $\v_1$ and $\v_2$: now
\be
\frac{t^v_2}{t^v_1}=\frac{r^v_2}{r^v_1}=\frac{h^2}{12}.
\ee

\item $R_1^A=0,\D$ orthotropic and $\delta_A=\delta_D=0$: in this case we get
\be
\begin{array}{c}
t^u=\dfrac{T_\gamma}{4T_1},\  r^u=0,\medskip\\ 
t^v_1=0,\ r^v_1=0,\medskip\\ 
t^v_2=0,\ r^v_2=\dfrac{h}{4}{R_0^B} {R_1^D} \left|{\dfrac{{T_\gamma }-(-1)^{\lambda } {T_1}\ \rho
}{\left(2 {R_1^D}^2-({R_0^D}+{T_0}) {T_1}\right) ({R_0^A}+{T_0})+3 {R_0^B}^2 {T_1} }}\right|,\medskip\\ 
\phi_2^v=0,\medskip\\
t^w=\dfrac{1}{4}\dfrac{ 3{R_0^B}^2T_\gamma-(R_0^A+T_0)-\left((R_0^D+T_0)T_\gamma-2(-1)^\lambda{R_1^D}^2\rho\right)     }{  \left(2 {R_1^D}^2-({R_0^D}+{T_0}) {T_1}\right) ({R_0^A}+{T_0})+3 {R_0^B}^2 {T_1}   },\medskip\\
r^w=\dfrac{1}{2}R_1^D(R_0^A+T_0)\left|\dfrac{T_\gamma-(-1)^\lambda T_1\ \rho}{\left(2 {R_1^D}^2-({R_0^D}+{T_0}) {T_1}\right) ({R_0^A}+{T_0})+3 {R_0^B}^2 {T_1}}\right|,\medskip\\
\phi^w=0.
\end{array}
\ee
So, $\u$ is isotropic, like $\U$, because $R_1^A=0\Rightarrow R^U=0$, and $\v_1=\gr{O}$: a change of temperature $t$ will produce a homogeneous state of in-plane forces and deformations but no curvature. Moreover, $\v_2$ is purely anisotropic, i.e. its isotropic part vanishes, $t^v_2=0$; as a consequence, a plate submitted to a through-the-thickness change of temperature $\nabla t$, by effect of the thermal coupling will stretch in all the directions but without change of surface. Also, because the isotropic part of $\w$ is not null, the mean curvature will be different from zero. This kind of laminates belongs to the set of {\it warp free thermally stable laminates}, cf. \cite{vannucci12joe1}: these are coupled laminates that do not warp as an effect of a change of temperature $t$, because $\v_1=\gr{O}$. 

\item $R_1^D=0,\A$ orthotropic and $\delta_A=\delta_D=0$:  this case is the dual of the previous one, obtained swapping $\A$ and $\D$, so we get
\be
\begin{array}{c}
t^u=\dfrac{1}{4}\dfrac{ 3{R_0^B}^2T_\gamma-(R_0^D+T_0)\left((R_0^A+T_0)T_\gamma-2(-1)^\lambda{R_1^A}^2\rho\right)     }{  \left(2 {R_1^A}^2-({R_0^A}+{T_0}) {T_1}\right) ({R_0^D}+{T_0})+3 {R_0^B}^2 {T_1}   },\medskip\\
r^u=\dfrac{1}{2}R_1^A(R_0^D+T_0)\left|\dfrac{T_\gamma-(-1)^\lambda T_1\ \rho}{\left(2 {R_1^A}^2-({R_0^A}+{T_0}) {T_1}\right) ({R_0^D}+{T_0})+3 {R_0^B}^2 {T_1}}\right|,\medskip\\
\phi^u=0,\medskip\\
t^v_1=0,\ r^v_1=\dfrac{h}{4}{R_0^B} {R_1^A} \left|{\dfrac{{T_\gamma }-(-1)^{\lambda } {T_1}\ \rho}{\left(2 {R_1^A}^2-({R_0^A}+{T_0}) {T_1}\right) ({R_0^D}+{T_0})+3 {R_0^B}^2 {T_1} }}\right|,\medskip\\ 
\phi_1^v=0,\medskip\\
t^v_2=0,\ r^v_2=0,\medskip\\ 
t^w=\dfrac{T_\gamma}{4T_1},\  r^w=0.
\end{array}
\ee
The remarks done for the previous case can be rephrased {\it verbatim}, changing  extension with bending. In particular, in this case a change of temperature $t$ will warp the plate, while a gradient $\nabla t$ will not stretch it: this is a case of {\it extension free thermally stable laminates}, see again \cite{vannucci12joe1}. For practical applications, if warping is sought for as a special effect, this case is very interesting.

\item $R_1^A=R_1^D,\delta_A=\delta_D=0$: this is the combination of the previous two cases and of course we get
\be
t^u=t^w=\dfrac{T_\gamma}{4T_1},\  r^u=r^w=0, \v_1=\v_2=\gr{O}:
\ee
this is a case of {\it warp- and extension-free thermally stable laminate}, cf. \cite{vannucci12joe1} once more. It is noticeable that now $\A\neq\D$, in general, but, cf. eq. (\ref{eq:formuleUVW}), $\U=\W$ and, by the last result above, $\u=\w$, and all of them are isotropic: the laminate is thermally quasi-homogeneous but not elastically. We remark that it is very simple to obtain such laminates: it is sufficient, but not necessary, to use square symmetric layers, i.e. reinforced by balanced fabrics. 

\item $R_0^A=R_0^D=0$: when the other special orthotropy is concerned, the results are not the same; in particular, in this case and also in all the two other ones, with only one among $R_0^A$ and $R_0^D$ that vanishes, the four compliance tensors $\u,\v_1,\v_2$ and $\w$ have always an isotropic and an anisotropic phase different from zero: no special situations exist. The results are:
\be
\begin{array}{c}
t^u={\dfrac{\left({T_0}^2-3 {R_0^B}^2\right) {T_1} {T_\gamma }+2 (-1)^{\lambda} {R_1^A}^2 \left(2 {R_1^D}^2-{T_0} {T_1}\right) \rho -2 {R_1^D}^2 {T_0} {T_\gamma }}{4\tau_0}},\medskip\\
r^u=\left|\dfrac{R_1^A (T_0 T_1-2 {R_1^D}^2)(T_\gamma-(-1)^\lambda T_1 \ \rho)}{2\tau_0}\right|,\medskip\\
\phi^u=0,
\end{array}
\ee
\be
\begin{array}{c}
t^v_1={-\dfrac{3}{h}R_0^BR_1^AR_1^D\dfrac{T_\gamma-(-1)^\lambda T_1 \ \rho}{\tau_0}},\medskip\\
r^v_1=\dfrac{3}{h}R_0^BR_1^AT_1\left|\dfrac{T_\gamma-(-1)^\lambda T_1 \ \rho}{\tau_0}\right|,\medskip\\
\phi^v_1=0,
\end{array}
\ee
\be
\begin{array}{c}
t^v_2={-\dfrac{h}{4}R_0^BR_1^AR_1^D\dfrac{T_\gamma-(-1)^\lambda T_1 \ \rho}{\tau_0}},\medskip\\
r^v_2=\dfrac{h}{4}R_0^BR_1^DT_1\left|\dfrac{T_\gamma-(-1)^\lambda T_1 \ \rho}{\tau_0}\right|,\medskip\\
\phi^v_2=0,
\end{array}
\ee
\be
\begin{array}{c}
t^w={\dfrac{\left({T_0}^2-3 {R_0^B}^2\right) {T_1} {T_\gamma }+2 (-1)^{\lambda} {R_1^D}^2 \left(2 {R_1^A}^2-{T_0} {T_1}\right) \rho -2 {R_1^A}^2 {T_0} {T_\gamma }}{4\tau_0}},\medskip\\
r^w=\left|\dfrac{R_1^D (T_0 T_1-2 {R_1^A}^2)(T_\gamma-(-1)^\lambda T_1 \ \rho)}{2\tau_0}\right|,\medskip\\
\phi^w=0,
\end{array}
\ee
with
\be
\tau_0=4 {R_1^A}^2 {R_1^D}^2-2 \left({R_1^A}^2+{R_1^D}^2\right) {T_0}{T_1}+ \left(-3 {R_0^B}^2+{T_0}^2\right) {T_1}^2.
\ee
To remark that the ratios $|r_1^v/t_1^v|, |r_2^v/t_2^v|,t_2^v/t_1^v$ and $r_2^v/r_1^v$ take exactly the same values of the case $\# i.$

\item $R_0^A=R_1^A=0$; this is the case of a laminate having an isotropic $\A$, but not necessarily an isotropic $\Ac$, because $\B\neq\mathbb{O}$, see \cite{vannucci23a}. Also in this case, we impose $\D$ to be orthotropic and $\delta_D=0$. The results are:
\be
\begin{array}{c}
t^u=\dfrac{T_\gamma}{4T_1},\  r^u=0,\medskip\\ 
t^v_1=0,\ r^v_1=0,\medskip\\ 
t^v_2=0,\ r^v_2=\dfrac{h}{4}{R_0^B} {R_1^D} \left|{\dfrac{{T_\gamma }-(-1)^{\lambda } {T_1}\ \rho
}{2 {R_1^D}^2 T_0+3 {R_0^B}^2 T_1-T_0 (R_0^D+T_0) T_1}}\right|,\medskip\\ 
\phi_2^v=0,\medskip\\
t^w=\dfrac{1}{4}\dfrac{ \left[3{R_0^B}^2-(R_0^D+T_0)T_0\right]T_\gamma+2(-1)^\lambda{R_1^D}^2T_0\ \rho     }{  2 {R_1^D}^2 T_0+3 {R_0^B}^2 T_1-T_0 (R_0^D+T_0) T_1 },\medskip\\
r^w=\dfrac{1}{2}R_1^DT_0\left|\dfrac{T_\gamma-(-1)^\lambda T_1\ \rho}{2 {R_1^D}^2 T_0+3 {R_0^B}^2 T_1-T_0 (R_0^D+T_0) T_1}\right|,\medskip\\
\phi^w=0.\\
\end{array}
\ee
The situation is very similar to that of case $\# iii.$ and the remarks are the same.\medskip
\end{enumerate}
To end this part, we notice that other combinations also exist, but less important for applications; one of them, the one with $\A=\D$ and isotropic, has the same result as the case $\# v.$; this should not be surprising, because all the thermoelastic tensors depend exclusively on the anisotropic phase varying with $2\theta$. So, substantially, adding to the condition that this phase vanishes, i.e. to case $\# v.$, another condition prescribing the nullity of the phase varying like $4\theta$ gives nothing more to the thermoelastic tensors, which is in accordance with the empirical {\it Neumann's principle}, \cite{love,nye}.

\subsection{Coupled laminates with $\V\neq\gr{O}$}
The case of $\V\neq\gr{O}$ is more complicate to be analyzed, but it is easier to find a stacking sequence, as they do not have to satisfy the condition $R_1^B=0$: in practice, all the coupled sequences are suitable. Reconsidering the cases of the previous Section, we can see a fundamental difference: when $\V\neq\gr{O}$, all the thermal tensors have an isotropic and an anisotropic phase, in all the cases but one: when $R_0^B=R_1^A=R_1^D=0$. That is why we present in the following only two cases: the same case $i.$ of the previous Section and the one with $R_0^B=R_1^A=R_1^D=0$.
\begin{enumerate}[i.]
\item $\A\neq\D$ but $k=0$ orthotropic and aligned with $\B\ (\delta_A=\delta_D=0)$: now we obtain
\be
\hspace{-3mm}
\begin{array}{c}
\begin{split}
t^u&=\dfrac{1}{4\psi}\left\{\left[12 R_0^B R_1^B R_1^D-6 {R_1^B}^2 (T_0+R_0^D)-3 {R_0^B}^2 T_1+R_0^D (R_0^A+T_0) T_1+\right.\right.\\
&\left.+(R_0^A+T_0) \left(T_0 T_1-2 {R_1^D}^2\right)\right] T_\gamma+2 (-1)^\lambda \left[2 \left(R_1^A R_1^D-3 {R_1^B}^2\right)^2-\right.\\
&\left.\left.-\left(R_0^D {R_1^A}^2-6 R_0^B R_1^A R_1^B+{R_1^A}^2 T_0+3 {R_1^B}^2 (R_0^A+T_0)\right) T_1\right] \rho\right\},
\end{split}\medskip\\
r^u=\left|\dfrac{\left[6 {R_1^B}^2 R_1^D-2 R_1^A {R_1^D}^2-3 R_0^B R_1^B T_1+R_1^A (R_0^D+T_0) T_1\right] (T_\gamma-(-1)^\lambda T_1 \rho)}{2\psi} \right|,\medskip\\
\phi^u=0,
\end{array}
\ee
\be
\hspace{-3mm}
\begin{array}{c}
\begin{split}t^v_1=\dfrac{3}{h}\dfrac{1}{\psi}&\left[R_0^D R_1^A R_1^B+R_0^A R_1^B R_1^D-R_0^B \left(3 {R_1^B}^2+R_1^A R_1^D\right)+\right.\\
&\left.+R_1^B (R_1^A+R_1^D) T_0\right] (T_\gamma-(-1)^\lambda T_1 \rho),\end{split}\medskip\\
r^v_1=\dfrac{3}{h}\left|\dfrac{\left[6 {R_1^B}^3+R_0^B R_1^A T_1-R_1^B (2 R_1^A R_1^D+(R_0^A+T_0) T_1)\right] (T_\gamma-(-1)^\lambda T_1 \rho)}{\psi}\right|,\medskip\\
\phi^v_1=0,
\end{array}
\ee
\be
\hspace{-3mm}
\begin{array}{c}
\begin{split}t^v_2=\dfrac{h}{4}\dfrac{1}{\psi}&\left[R_0^D R_1^A R_1^B+R_0^A R_1^B R_1^D-R_0^B \left(3 {R_1^B}^2+R_1^A R_1^D\right)+\right.\\
&\left.+R_1^B (R_1^A+R_1^D) T_0\right] (T_\gamma-(-1)^\lambda T_1 \rho),\end{split}\medskip\\
r^v_2=\dfrac{h}{4}\left|\dfrac{\left[ 6 {R_1^B}^3+R_0^B R_1^D T_1-R_1^B \left(2 R_1^A R_1^D+(R_0^D+T_0) T_1\right)\right] (T_\gamma-(-1)^\lambda T_1 \rho)}{\psi}\right|,\medskip\\
\phi^v_2=0,
\end{array}
\ee
\be
\hspace{-3mm}
\begin{array}{c}
\begin{split}
t^w&=\dfrac{1}{4\psi}\left\{\left[12 R_0^B R_1^B R_1^A-6 {R_1^B}^2 (T_0+R_0^A)-3 {R_0^B}^2 T_1+R_0^A (R_0^D+T_0) T_1+\right.\right.\\
&\left.+(R_0^D+T_0) \left(T_0 T_1-2 {R_1^A}^2\right)\right] T_\gamma+2 (-1)^\lambda \left[2 \left(R_1^D R_1^A-3 {R_1^B}^2\right)^2-\right.\\
&\left.\left.-\left(R_0^A{R_1^D}^2-6 R_0^B R_1^D R_1^B+{R_1^D}^2 T_0+3 {R_1^B}^2 (R_0^D+T_0)\right) T_1\right] \rho\right\},
\end{split}\medskip\\
r^w=\left|\dfrac{\left[6 {R_1^B}^2 R_1^A-2 R_1^D{R_1^A}^2-3 R_0^B R_1^B T_1+R_1^D (R_0^A+T_0) T_1\right] (T_\gamma-(-1)^\lambda T_1 \rho)}{2\psi} \right|,\medskip\\
\phi^w=0,
\end{array}
\ee
with
\be
\begin{split}
\psi&=36 {R_1^B}^4+12 R_0^B R_1^B (R_1^A+R_1^D) T_1+2 {R_1^A}^2 \left[2 {R_1^D}^2-(R_0^D+T_0) T_1\right]-\\
&-6 {R_1^B}^2 \left[4 R_1^A R_1^D+(R_0^A+R_0^D+2 T_0) T_1\right]+\\
&+T_1 \left[(R_0^A+T_0)\left ((R_0^D+T_0) T_1-2 {R_1^D}^2\right)-3 {R_0^B}^2 T_1\right].
\end{split}
\ee
This is the most general case that we consider; it is worth noting that once more eq. (\ref{eq:ratiotv}) is valid: the value $h^2/12$ for the ratio of the two isotropic parts of $\v_2$ and $\v_1$ has hence always this value, unless one of the two tensors is purely anisotropic, i.e. when $t^v_j=0,j\in\{1,2\}$. The ratio $r^v_2/r^v_1$ takes a value different from that of eq. (\ref{eq:ratiorv}), because $R_1^B\neq0$; however, for TQHCL, i.e. if $R_0^A=R_0^D,R_1^A=R_1^D$, then it takes again the value $h^2/12$.  The case of TQHCL plates can be easily obtained putting in the above equations exactly $R_0^A=R_0^D,R_1^A=R_1^D$.

\item Laminates with $R_1^A=R_1^D=R_0^B=0$: this is the minimal requirement to get tensors with an elastic phase that vanishes; in particular, it is 
\be
\U=\W=\left\{\begin{array}{c}T_\gamma\\T_\gamma\\0\end{array}\right\},\ \ \V=\left\{\begin{array}{c}(-1)^\lambda R_1^B\ \rho\\(-1)^\lambda R_1^B\ \rho\\0\end{array}\right\}
\ee
and
\be
\begin{split}
&t^u=\frac{(R_0^D+T_0) T_\gamma-6 (-1)^\lambda {R_1^B}^2 \rho}{4\left[(R_0^D+T_0) T_1-6 {R_1^B}^2\right]},\ \ r^u=0,\\
&t^v_1=0,\ \ r^v_1=\frac{3}{h}\left|\frac{R_1^B (T_\gamma-(-1)^\lambda T_1\ \rho)}{(R_0^D+T_0) T_1-6 {R_1^B}^2}\right|,\ \ \phi^v_1=0,\\
&t^v_2=0,\ \ r^v_2=\frac{h}{4}\left|\frac{R_1^B (T_\gamma-(-1)^\lambda T_1\ \rho)}{(R_0^A+T_0) T_1-6 {R_1^B}^2}\right|,\ \ \phi^v_2=0,\\
&t^w=\frac{(R_0^A+T_0) T_\gamma-6 (-1)^\lambda {R_1^B}^2 \rho}{4\left[(R_0^A+T_0) T_1-6 {R_1^B}^2\right]},\ \ r^w=0.
\end{split}
\ee
We see hence that $\u$ and $\w$ are isotropic, while $\v_1$ and $\v_2$ purely anisotropic ($r_1^v=r_2^v=0$). These laminates are particularly interesting, because to obtain them is rather simple, it is sufficient to use balanced cross-ply sequences, i.e. laminates having an equal number of plies at $0^\circ$ and $90^\circ$, with a stacking sequence antisymmetric or belonging to the set of {\it quasi-trivial} laminates with $\C=\mathbb{O}$ and $\B\neq\mathbb{O}$, cf. \cite{vannucci01ijss,CompScTech01} or \cite{vannucci_libro}. 
\end{enumerate}

\section{Some examples}
We consider first the effects of a change of temperature $t\neq0,\ \nabla t=0$ on a coupled laminate. This is the most interesting situation for applications, more than the case of a $\nabla t\neq0$, because normally anisotropic laminates are thin plates so a change of temperature rapidly is uniform through the thickness. Typically, this is the practical situation of all the plates fabricated through a process of curing pre-preg layers: once fabricated, the laminate is cooled down to the room temperature. Let us imagine such a laminate free to deform, not acted upon by external forces, so that all the internal actions, in-plane forces and bending/twisting moments, and their deformation effects are absent. With these assumptions eq. (\ref{eq:fundlawinversetherm}) gives
\be
\beps=t\ \u,\ \ \bk=t\ \v_1,
\ee
i.e. it is stretched and bent by $t$. We consider the following 12-ply laminates:
\be
\label{eq:sequenze}
\begin{split}
&[ 0^\circ\   0^\circ\   0^\circ\   0^\circ\   0^\circ\   0^\circ\   90^\circ\   90^\circ\   90^\circ\   90^\circ\   90^\circ\   90^\circ],\\
&[ 0^\circ\   90^\circ\   90^\circ\   0^\circ\   0^\circ\   90^\circ\   0^\circ\   90^\circ\   90^\circ\   0^\circ\   0^\circ\   90^\circ],\\
&[ 0^\circ\   0^\circ\   90^\circ\   90^\circ\   0^\circ\   90^\circ\   0^\circ\   90^\circ\   0^\circ\   0^\circ\   90^\circ\   90^\circ],\\
&[ 0^\circ\   90^\circ\   90^\circ\   0^\circ\   0^\circ\   0^\circ\   90^\circ\   0^\circ\   90^\circ\   90^\circ\   90^\circ\   0^\circ],\\
&[ 0^\circ\   90^\circ\   90^\circ\   0^\circ\   0^\circ\   0^\circ\   90^\circ\   0^\circ\   90^\circ\   90^\circ\   90^\circ\   0^\circ],\\
&[ 0^\circ\   90^\circ\   90^\circ\   0^\circ\   0^\circ\   90^\circ\   0^\circ\   0^\circ\   90^\circ\   90^\circ\   90^\circ\   0^\circ].
\end{split}
\ee
All of them are TQHCL, with six layers at $0^\circ$ and the six others at $90^\circ$; the first three sequences are antisymmetric, the last three completely asymmetric. Tensors $\A,\D,\U$ and $\W$ are the same for all the laminates, but not $\B$ nor $\V$ and, by consequence, also $\Ac,\Bc,\Dc,\u,\v_1,\v_2$ and $\w$ change for each stacking sequence. Moreover, $\A=\D,\U=\W,\u=\w$ and $R_1^A=R_1^D=R_0^B=0$ for all the cases.

The material for the plies is carbon-epoxy  T300/5208, \cite{TsaiHahn}, whose characteristics are: 
\begin{center}
$E_1= 181000$ MPa, $E_2=10300$ MPa, $G_{12}=7170$ MPa, $\nu_{12}=0.28$,\medskip\\ $\alpha_1=2\times10^{-6}\ ^\circ$C$^{-1}$, $\alpha_2=2.25\times10^{-3}\ ^\circ$C$^{-1}$, \medskip
\end{center}
which gives \medskip
\begin{center}
$T_0=26880$ MPa, $T_1=24740$ MPa, $R_0=19710$ MPa, $R_1=21430$ MPa, $\Phi_0=\Phi_1=0$,\medskip \\$T_\gamma=15.1\ ^\circ$C$^{-1}$MPa, $R_\gamma=8.21\ ^\circ$C$^{-1}$MPa, $\Phi_\gamma=0$, $\rho=3.83\hspace{-1mm}\times\hspace{-1mm}10^{-4}\ ^\circ$C$^{-1}$.\medskip
\end{center}

\begin{figure}[h]
\center
\includegraphics[width=.6\textwidth]{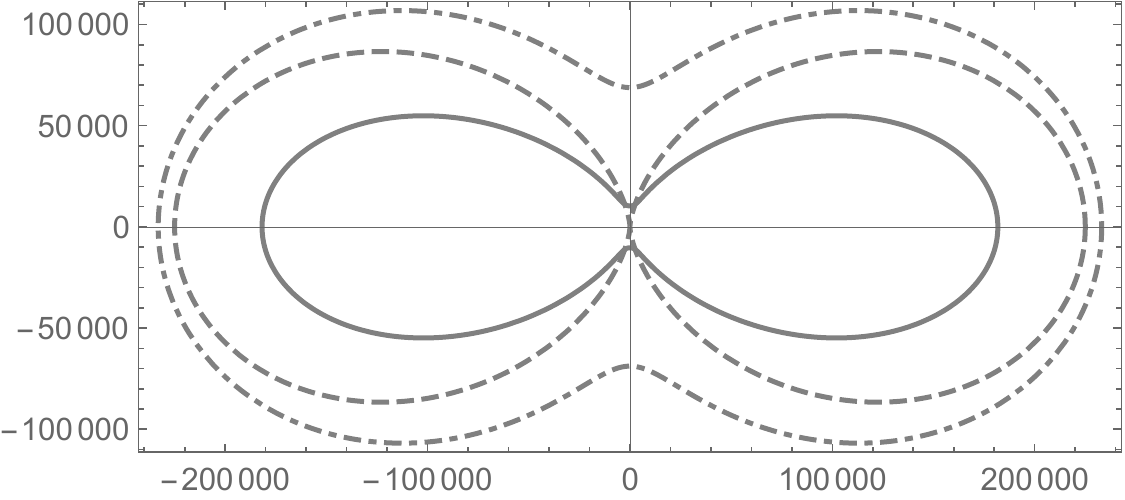}
\caption{Polar diagrams of $\Q_{11}(\theta)$, continuous line, $\alpha_1(\theta)\times10^8$, dashed line, $\gamma_1(\theta)\times10^4$, dotdashed line.}
\label{fig:4}
\end{figure}

\begin{figure}[h]
\center
\includegraphics[width=.6\textwidth]{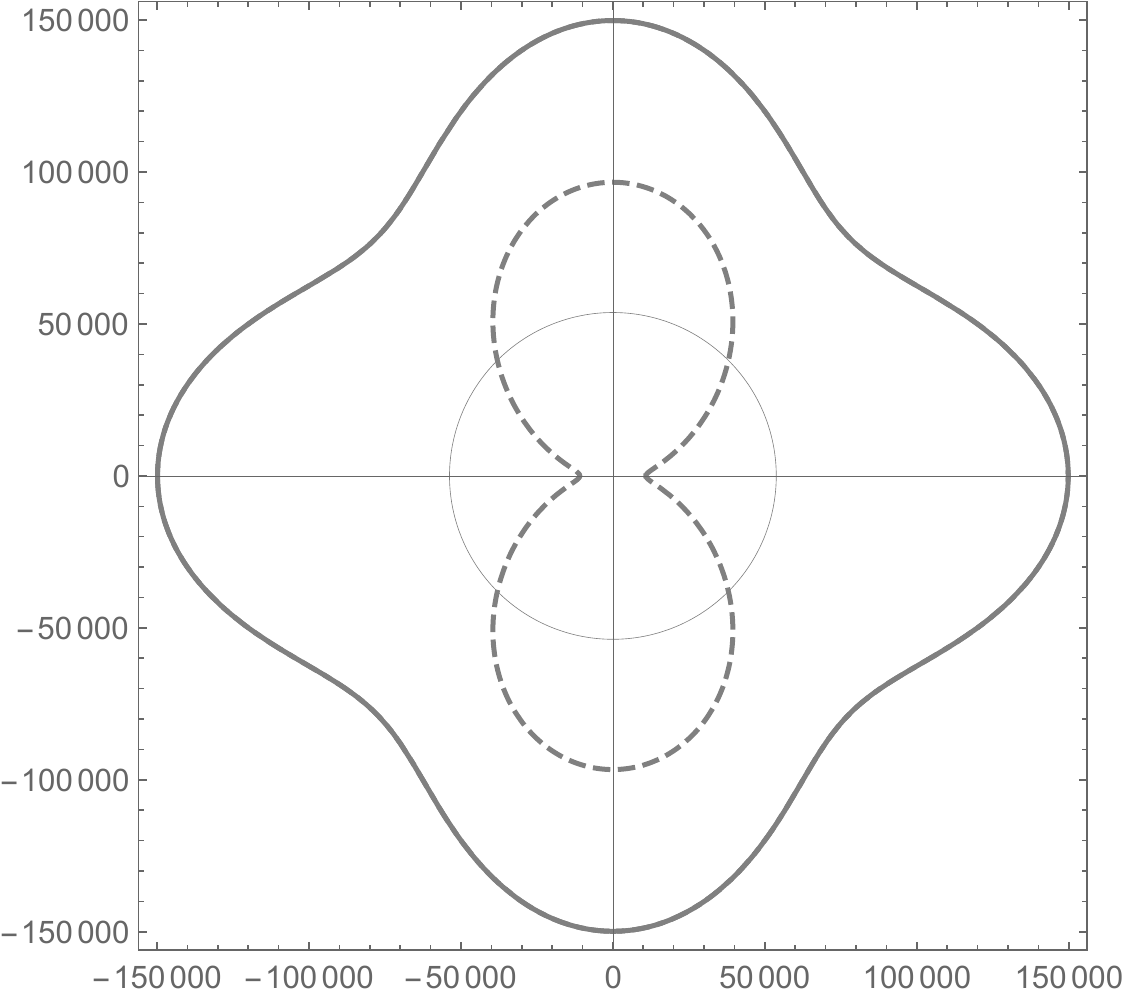}
\caption{Polar diagrams of $\A_{11}(\theta)=\D_{11}(\theta)$, continuous line, $\B_{11}(\theta)$, dashed line for the first laminate; the thin line is the zero-line: inside it, a diagram represents a negative value.}
\label{fig:5}
\end{figure}

\begin{figure}[h]
\center
\includegraphics[width=.6\textwidth]{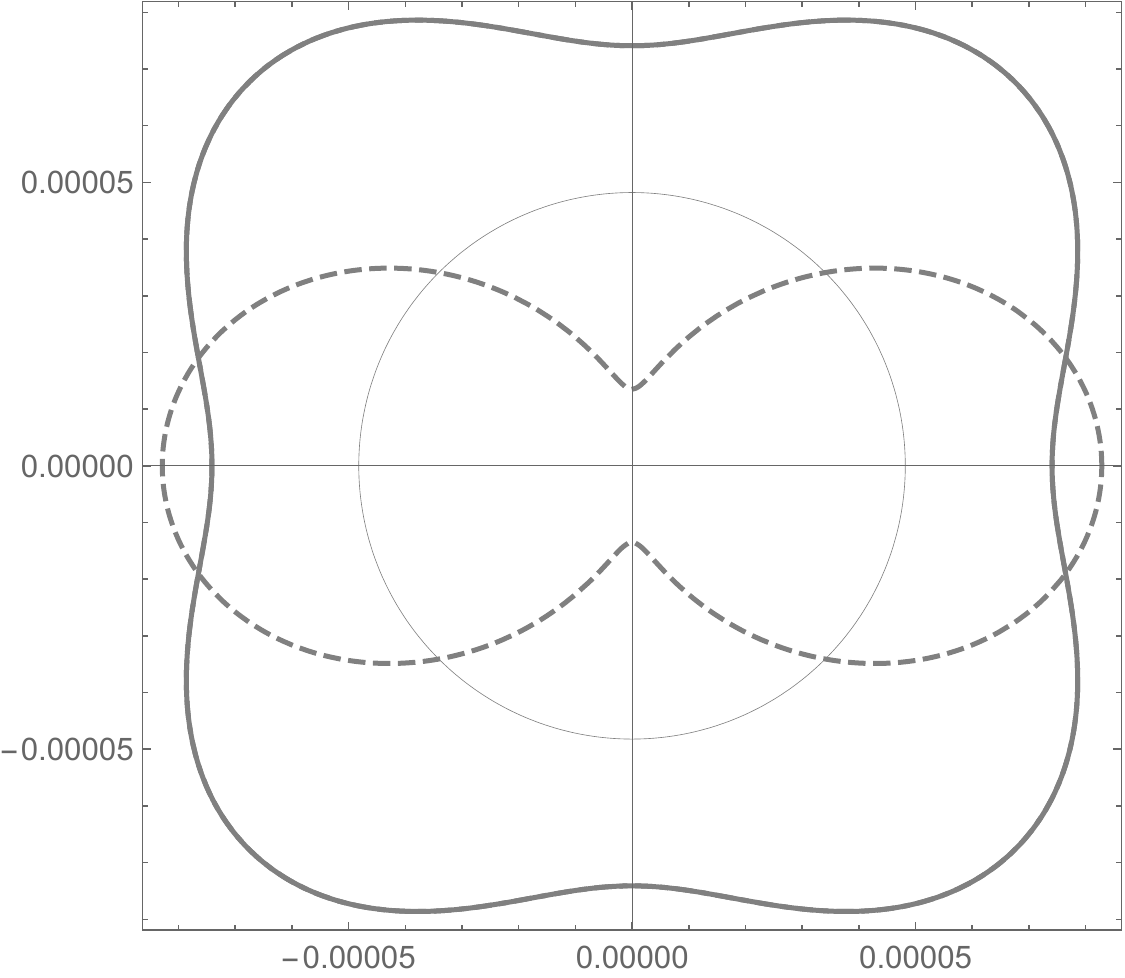}
\caption{Polar diagrams of $\Ac_{11}(\theta)=\Dc_{11}(\theta)$, continuous line, $\Bc_{11}(\theta)$, dashed line for the first laminate; the thin line is the zero-line: inside it, a diagram represents a negative value.}
\label{fig:6}
\end{figure}

\begin{figure}[h]
\center
\includegraphics[width=.6\textwidth]{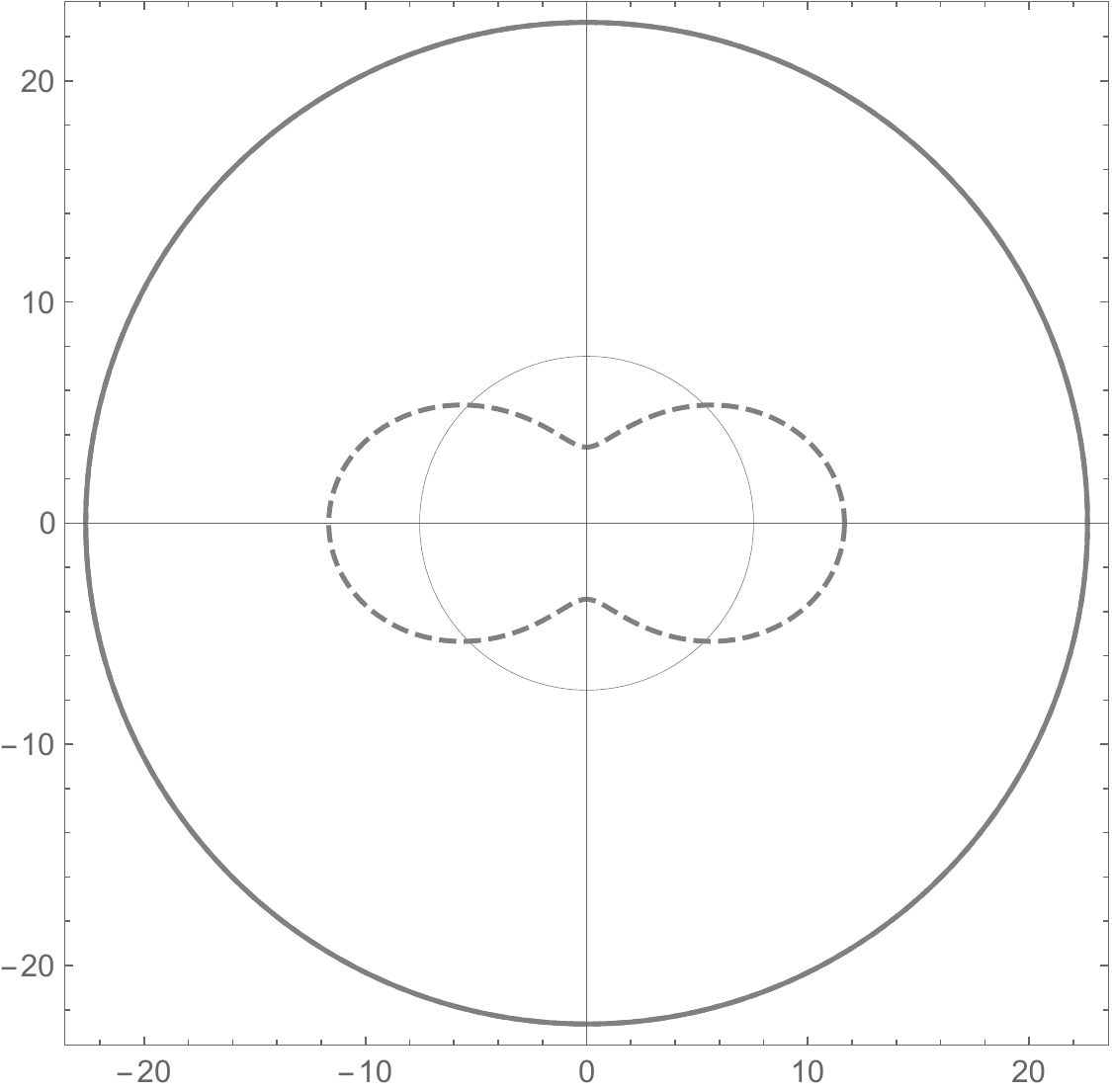}
\caption{Polar diagrams of $\U_{1}(\theta)=\W_{1}(\theta)$, continuous line, $\V_{1}(\theta)$, dashed line for the first laminate; the thin line is the zero-line: inside it, a diagram represents a negative value.}
\label{fig:7}
\end{figure}

\begin{figure}[h]
\center
\includegraphics[width=.6\textwidth]{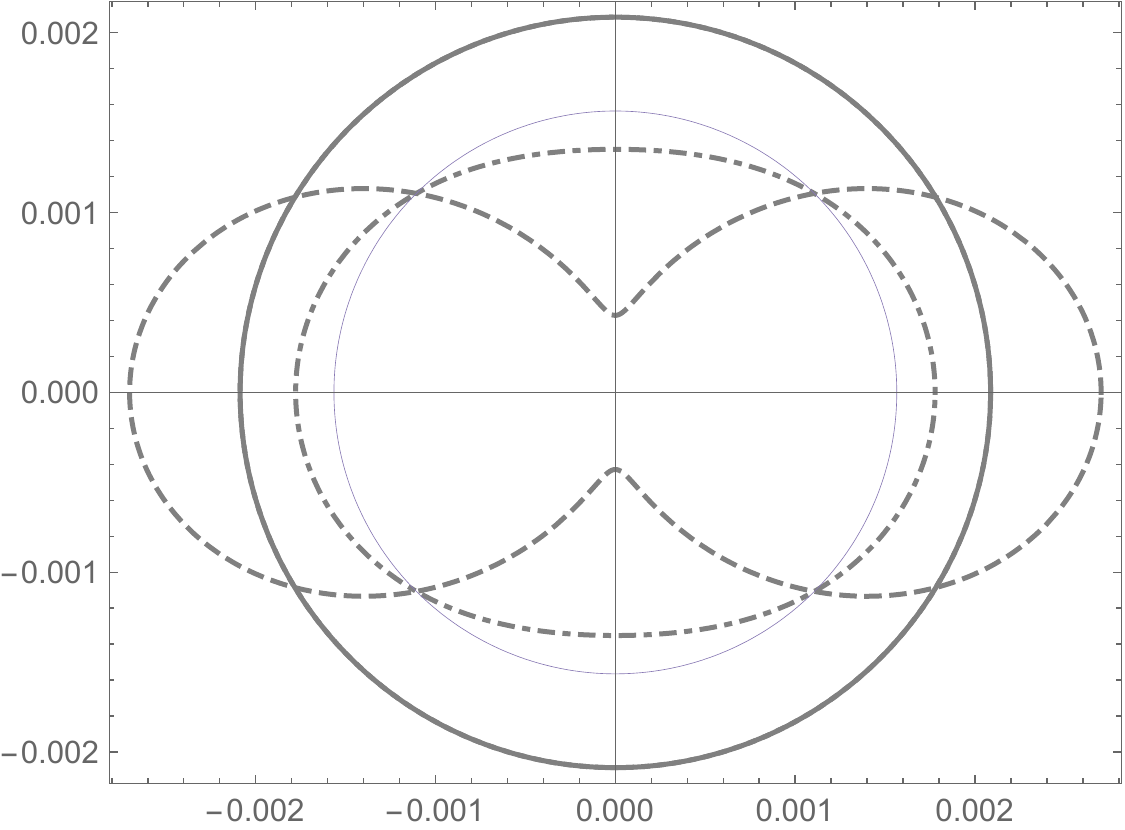}
\caption{Polar diagrams of $\u_{1}(\theta)=\w_{1}(\theta)$, continuous line, $(\v_{1})_1(\theta)$, dashed line, $(\v_{2})_1(\theta)$, dotdashed line for the first laminate; the thin line is the zero-line: inside it, a diagram represents a negative value.}
\label{fig:8}
\end{figure}

\begin{figure}[h]
\center
\includegraphics[width=.6\textwidth]{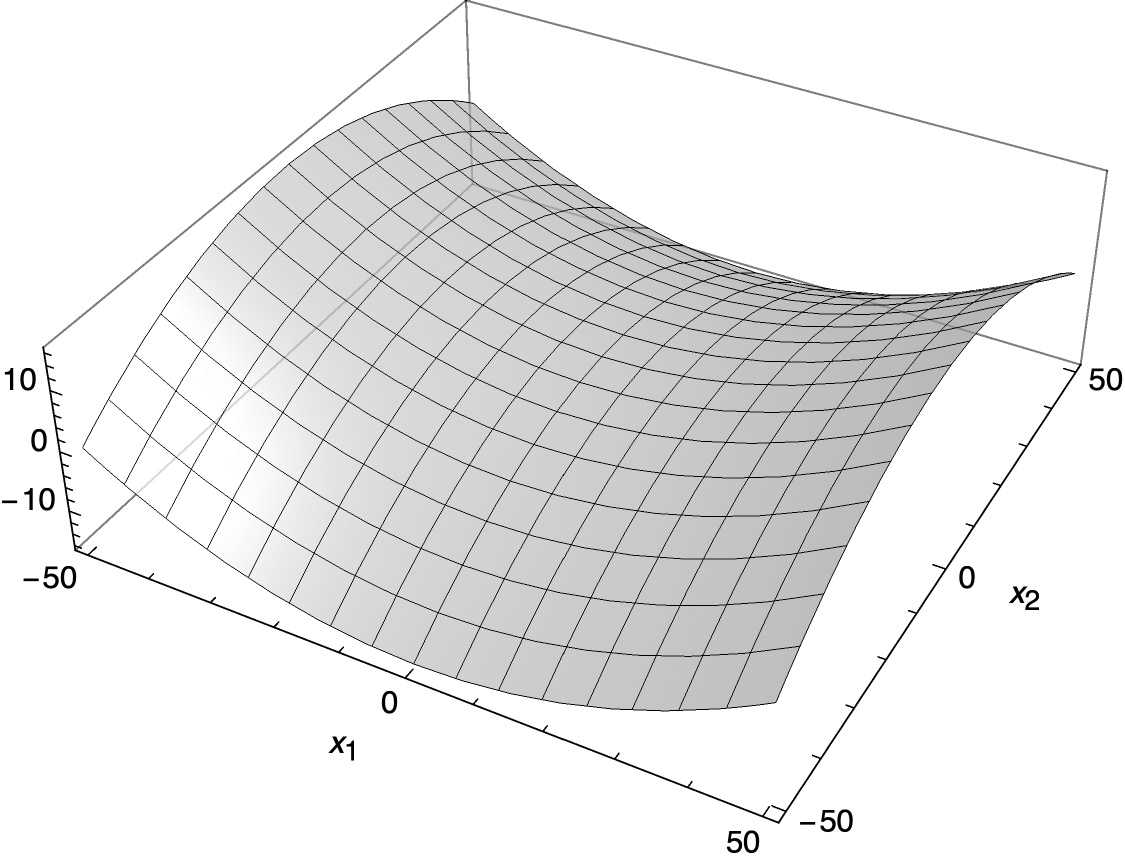}
\caption{Bending effect of a temperature change $t=50^\circ$C due to coupling $\v_1$ on a square plate with a side of 100 mm.}
\label{fig:9}
\end{figure}

\begin{figure}[h]
\center
\includegraphics[width=.6\textwidth]{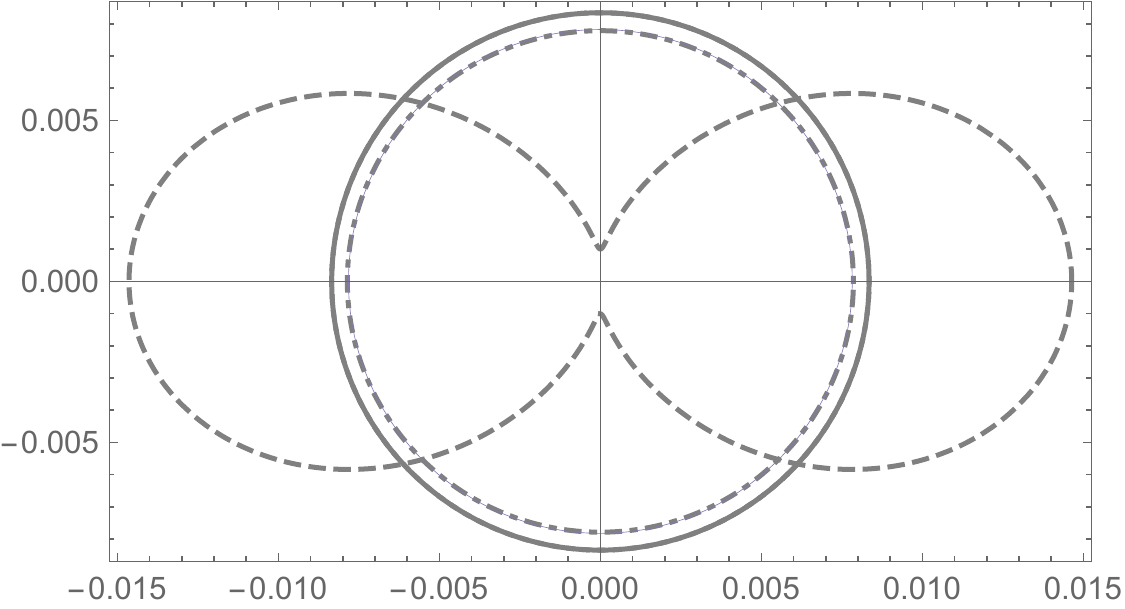}
\caption{Polar diagrams of $\u_{1}(\theta)=\w_{1}(\theta)$, continuous line, $(\v_{1})_1(\theta)$, dashed line, $(\v_{2})_1(\theta)$, dotdashed line, for the two-plies laminate; the thin line is the zero-line: inside it, a diagram represents a negative value.}
\label{fig:10}
\end{figure}

Let us consider in detail just two of the previous laminates. \medskip\\
{\bf First laminate}:  this is the first sequence in eq. (\ref{eq:sequenze}). It is  a completely antisymmetric sequence and  has the highest degree of coupling, as defined in \cite{vannucci01joe}. For this plate, we get:\medskip
\begin{itemize}
\item $\A=\D=\left[\begin{array}{ccc}9.61&0.29&0\\0.29&9.61&0\\0&0&1.43\end{array}\right]\times10^4$ MPa, $\B=\left[\begin{array}{ccc}-4.29&0&0\\0&4.29&0\\0&0&0\end{array}\right]\times10^4$ MPa,\\
$R_0^A=R_0^D=1.97\times10^4$ MPa, $R_1^B=1.07\times10^4$ MPa;\bigskip
\item $\Ac\hspace{-1mm}=\hspace{-1mm}\Dc\hspace{-1mm}=\hspace{-1mm}\left[\begin{array}{ccc}2.59&-0.08&0\\-0.08&2.59&0\\0&0&6.97\end{array}\right]\hspace{-.5mm}\times\hspace{-0.5mm}10^{-5}$ MPa$^{-1}$, $\Bc\hspace{-1mm}=\hspace{-1mm}\left[\begin{array}{ccc}3.47&0&0\\0&-3.47&0\\0&0&0\end{array}\right]\hspace{-0.5mm}\times\hspace{-0.5mm}10^{-5}$ MPa$^{-1}$,\\
$t_0^A=t_0^D=2.41\times10^{-5}$ MPa$^{-1}$, $t_1^A=t_1^D=6.28\times10^{-6}$ MPa$^{-1}$,\\
$r_0^A=r_0^D=1.08\times10^{-5}$ MPa$^{-1}$, $r_1^A=r_1^D=0, \phi_0^A=\phi_0^D=45^\circ,\phi_1^A=\phi_1^D=0$;\bigskip
\item$\U=\W=\left\{\begin{array}{c}15.1\\15.1\\0   \end{array}\right\}\ ^\circ$C$^{-1}$MPa, $\V=\left\{\begin{array}{c}4.11\\-4.11\\0   \end{array}\right\}\ ^\circ$C$^{-1}$MPa,\\
$T^U=T^W=15.1\ ^\circ$C$^{-1}$MPa, $R^U=R^W=0$, $T^V=0,R^V=4.11\ ^\circ$C$^{-1}$MPa; \bigskip
\item $\u=\w=\left\{\begin{array}{c}5.21\\5.21\\0   \end{array}\right\}\times10^{-4}\ ^\circ$C$^{-1}$, $\v_1=\left\{\begin{array}{c}1.14\\-1.14\\0   \end{array}\right\}\times10^{-3}\ ( ^\circ$C mm)$^{-1}$,\\ $\v_2=\left\{\begin{array}{c}2.13\\-2.13\\0   \end{array}\right\}\times10^{-4}\ ^\circ$C$^{-1}$mm,\\
$t^u=t^w=5.21\times10^{-4}\ ^\circ$C$^{-1}$, $r^u=r^w=0, t^v_1=t^v_2=0,$\\
$r^v_1=1.14\times10^{-3}\ ( ^\circ$C mm)$^{-1}$, $r^v_2=2.13\times10^{-4}\ ^\circ$C$^{-1}$mm.\medskip
\end{itemize}
The polar diagrams of the first component of the above tensors are shown in Figs. \ref{fig:4} to \ref{fig:8}.
The temperature change $t$ has two effects:
\begin{enumerate}[a)]
\item it causes an {\it extension}, i.e. an isotropic dilatation, if $t>0$, or contraction, if $t<0$, because $\u$ is isotropic,  regardless the fact that the plate is, elastically, anisotropic (actually, it is $R_0$-orthotropic). The amount of the extension  is
\be
e=t\ t^u=\frac{(R_0^D+T_0) T_\gamma-6 (-1)^\lambda {R_1^B}^2 \rho}{4\left[(R_0^D+T_0) T_1-6 {R_1^B}^2\right]}t,
\ee
for the case at hand $e=5.21\times10^{-4}\ t$;
\item by tensor coupling $\v_1$, $t$ bends the plate: the curvature tensor is purely anisotropic, like $\v_1$:
\be
\label{eq:curvature}
\bk=t\ \v_1=\frac{3}{h}\left|\frac{R_1^B (T_\gamma-(-1)^\lambda T_1\ \rho)}{(R_0^D+T_0) T_1-6 {R_1^B}^2}\right|\left\{\begin{array}{c}t\\-t\\0   \end{array}\right\}=\left\{\begin{array}{c}1.14\\-1.14\\0   \end{array}\right\}t\times10^{-3}\  \mathrm{ mm}^{-1}.
\ee
\end{enumerate}
The curvature effect of $t$ due to coupling is much more important than the extension one; curvature can be very important and the temperature change $t$ and the coupling tensor $\v_1$ give rise to a {\it bistable phenomenon}: such a plate can be in an equilibrium configuration also with a change of the sign of the two curvatures (this is due to the square-symmetry of the plate).  Moreover, the deformed surface of the plate is made of hyperbolic points, because the Gaussian curvature $K=\kappa_1\kappa_2$, is negative everywhere, and it is a minimal surface, because the mean curvature of the plate, $H=\frac{1}{2}(\kappa_1+\kappa_2)$, is null everywhere, \cite{pressley,vannucci_alg}. In Fig. \ref{fig:9} a sketch of the bent surface for a square plate with a side of 100 mm and $t=50^\circ$C.
If the stacking sequence is the same but the number of plies changes,  only   tensors $\v_1$ and $\v_2$ change, while all the other ones remain unchanged. In particular, if we reduce the plate to only two layers, the first one at $0^\circ$, the other one at $90^\circ$, we get\medskip
\begin{center}
$\v_1=\left\{\begin{array}{c}6.82\\-6.82\\0   \end{array}\right\}\times10^{-3}\ ( ^\circ$C mm)$^{-1}$,\ $\v_2=\left\{\begin{array}{c}3.55\\-3.55\\0   \end{array}\right\}\times10^{-5}\ ^\circ$C$^{-1}$mm,\medskip\\
$r^v_1=6.82\times10^{-3}\ ( ^\circ$C mm)$^{-1}$, $r^v_2=3.55\times10^{-5}\ ^\circ$C$^{-1}$mm,\medskip
\end{center}
so $\v_1$ increases of  six times, while $\v_2$ decreases of the same amount (the polar diagrams of the new tensors are shown in Fig. \ref{fig:10}). 
 However,  none of these quantities affects the value of the curvatures, eq. (\ref{eq:curvature}), but the plate's thickness, that decreases to $h=0.25$ mm, so the curvatures increases of six times; in Fig. \ref{fig:11} it is possible to appreciate the effect of the plate's thickness decrease. Actually, the effect of $t$ is inversely proportional to $h$, or, which is the same, to the plis number $n$: the less the thickness, the highest the curvature.
\begin{figure}[h]
\center
\includegraphics[width=.4\textwidth]{11.pdf}
\caption{Bending effect of a temperature change $t=50^\circ$C due to coupling $\v_1$ on a square plate with a side of 100 mm, for the two cases of a laminate with 12 (dark grey) or 2 plies (light grey).}
\label{fig:11}
\end{figure}

Let us now consider the effect of a temperature gradient, i.e. the case $t=0,\nabla t\neq0$; now, eq. (\ref{eq:fundlawinversetherm}) gives
\be
\beps=\nabla t\ \v_2,\ \ \kappa=\nabla t\ \w.
\ee
In particular, because $\v_2$ is purely anisotropic ($t_2^v=0$), we get
\be
\beps=\frac{h}{4}\left|\frac{R_1^B (T_\gamma-(-1)^\lambda T_1\ \rho)}{(R_0^A+T_0) T_1-6 {R_1^B}^2}\right|\left\{\begin{array}{c}\nabla t\\-\nabla t\\0   \end{array}\right\}=\left\{\begin{array}{c}2.13\\-2.13\\0   \end{array}\right\}\nabla t\times10^{-4}.
\ee
The in-plane deformation is hence isochoric  and homogeneous, because $\tr\beps=0$: the plate deforms in its plane without changing of surface.

\begin{figure}[h]
\center
\includegraphics[width=.5\textwidth]{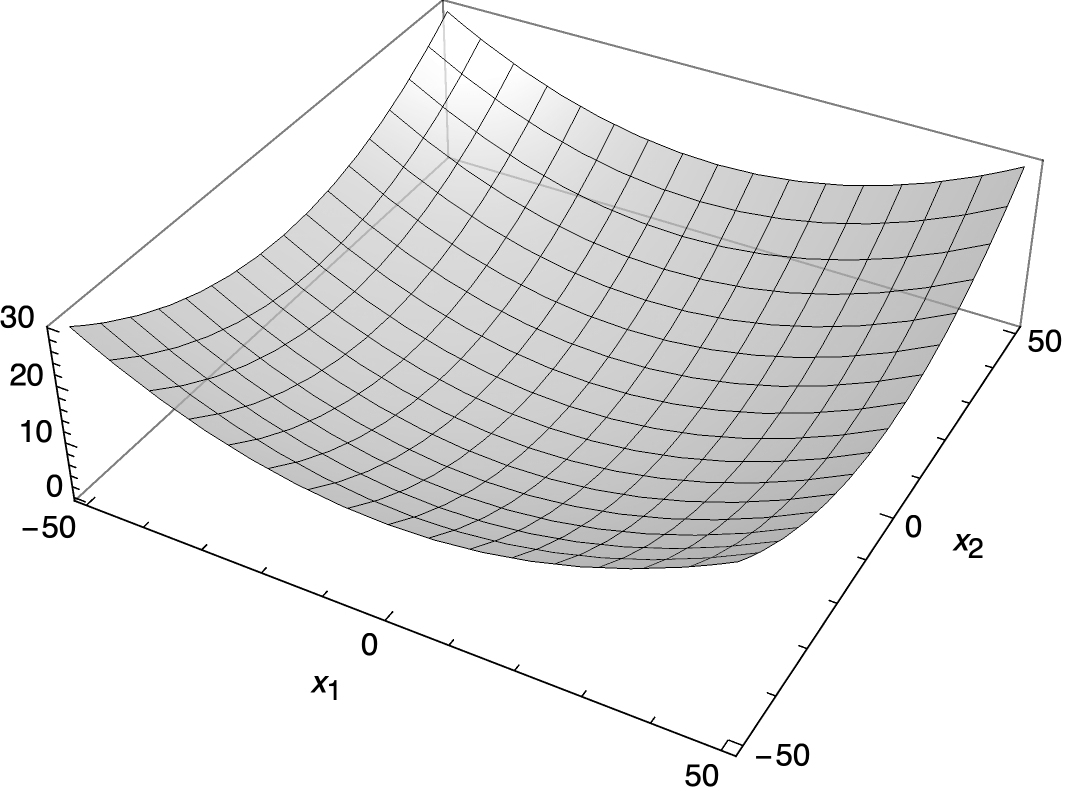}
\caption{Bending effect of a gradient of temperature  $\nabla t=50^\circ$C/mm due to  $\w$ on a square plate with a side of 100 mm.}
\label{fig:12}
\end{figure}

\begin{figure}[h]
\center
\includegraphics[width=.6\textwidth]{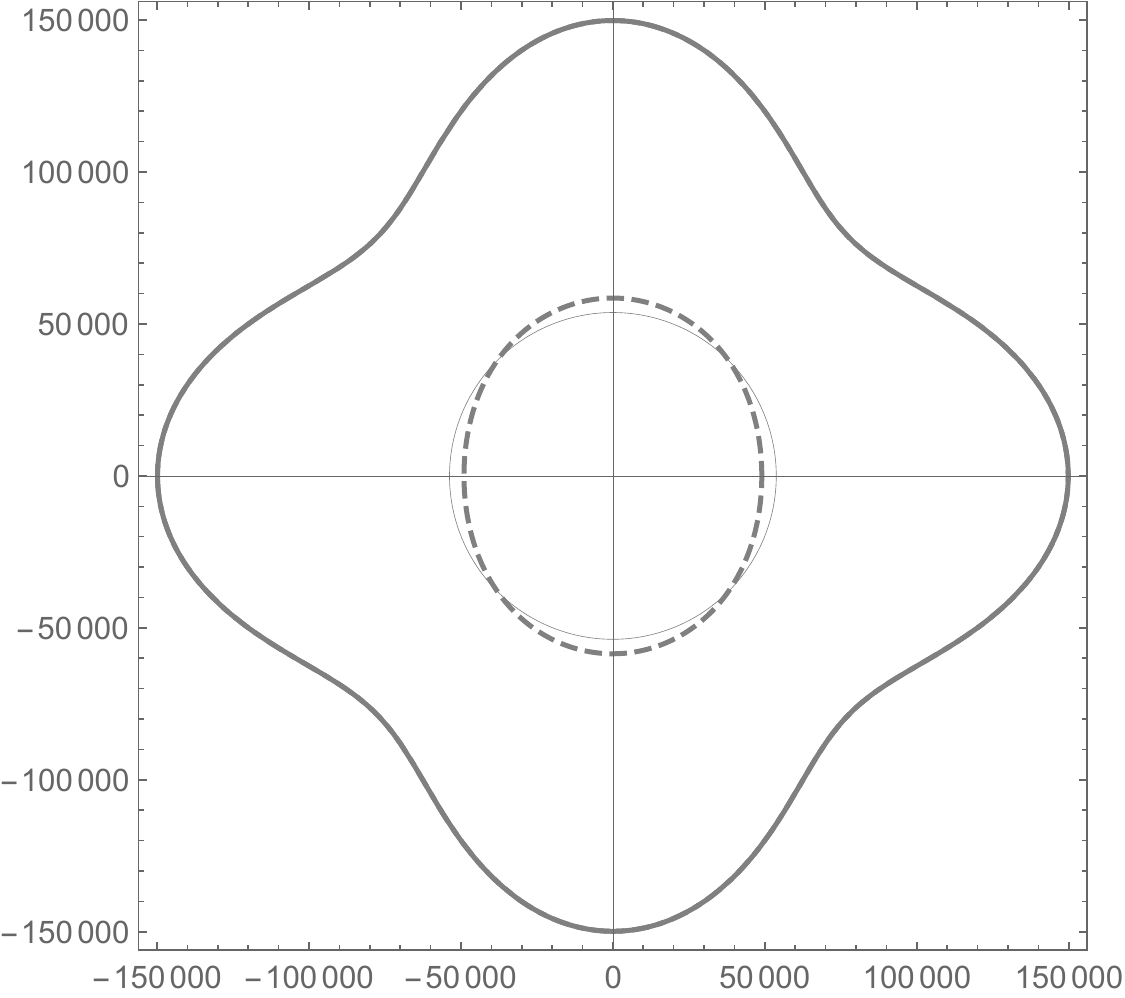}
\caption{Polar diagrams of $\A_{11}(\theta)=\D_{11}(\theta)$, continuous line, $\B_{11}(\theta)$, dashed line for the second laminate; the thin line is the zero-line: inside it, a diagram represents a negative value.}
\label{fig:13}
\end{figure}

\begin{figure}[h]
\center
\includegraphics[width=.6\textwidth]{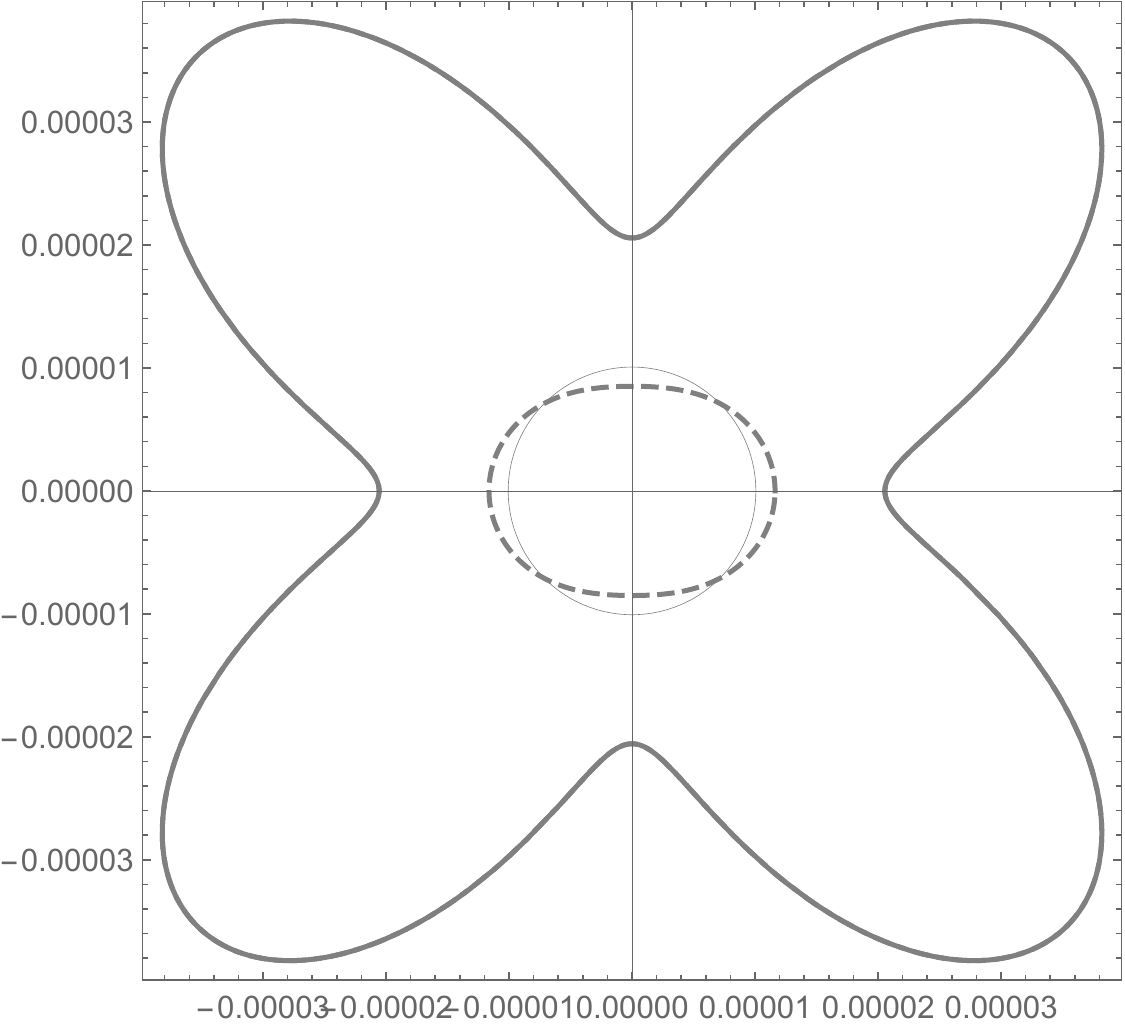}
\caption{Polar diagrams of $\Ac_{11}(\theta)=\Dc_{11}(\theta)$, continuous line, $\Bc_{11}(\theta)$, dashed line for the second laminate; the thin line is the zero-line: inside it, a diagram represents a negative value.}
\label{fig:14}
\end{figure}

\begin{figure}[h]
\center
\includegraphics[width=.6\textwidth]{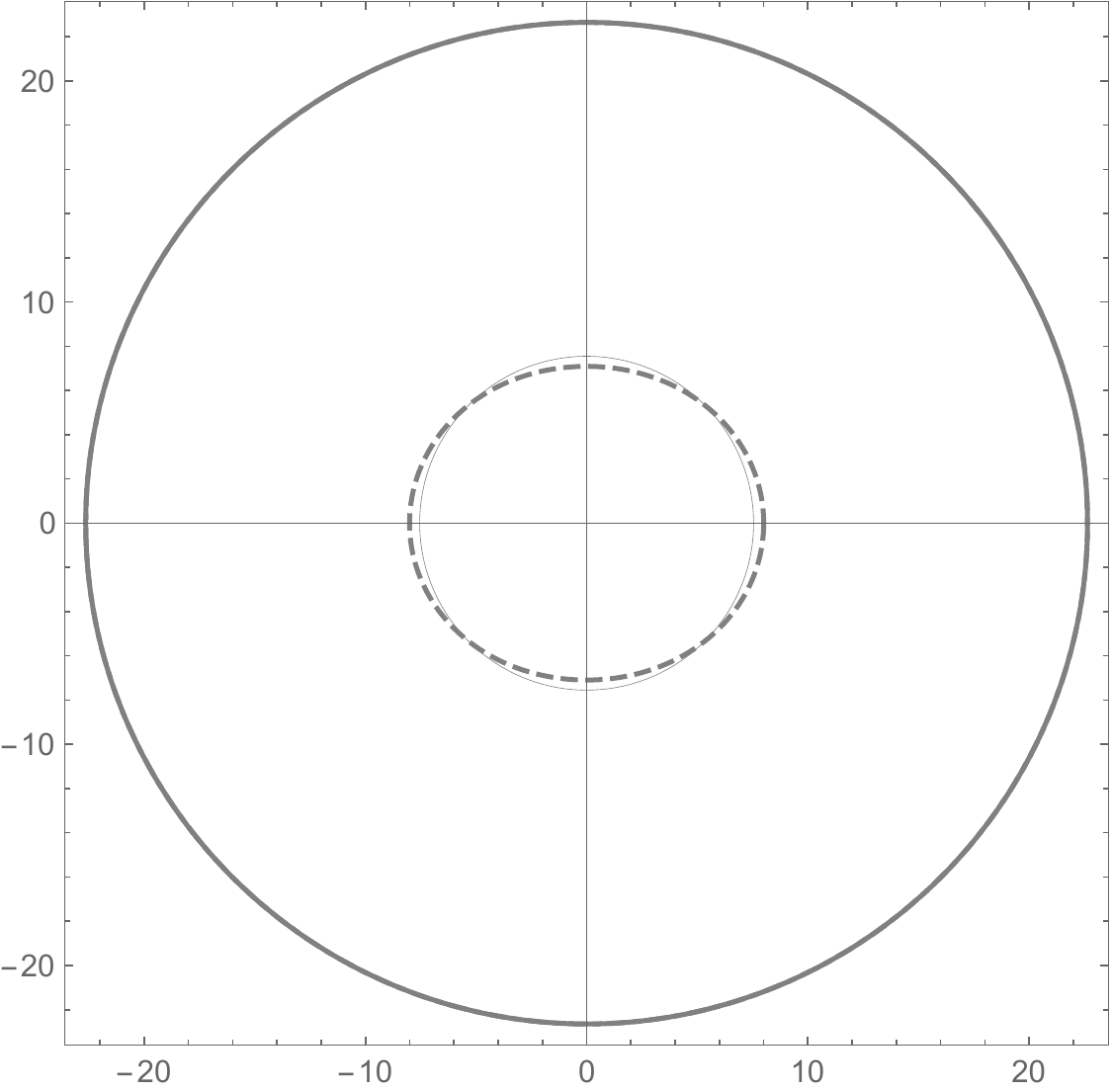}
\caption{Polar diagrams of $\U_{1}(\theta)=\W_{1}(\theta)$, continuous line, $\V_{1}(\theta)$, dashed line for the second laminate; the thin line is the zero-line: inside it, a diagram represents a negative value.}
\label{fig:15}
\end{figure}

\begin{figure}[h]
\center
\includegraphics[width=.6\textwidth]{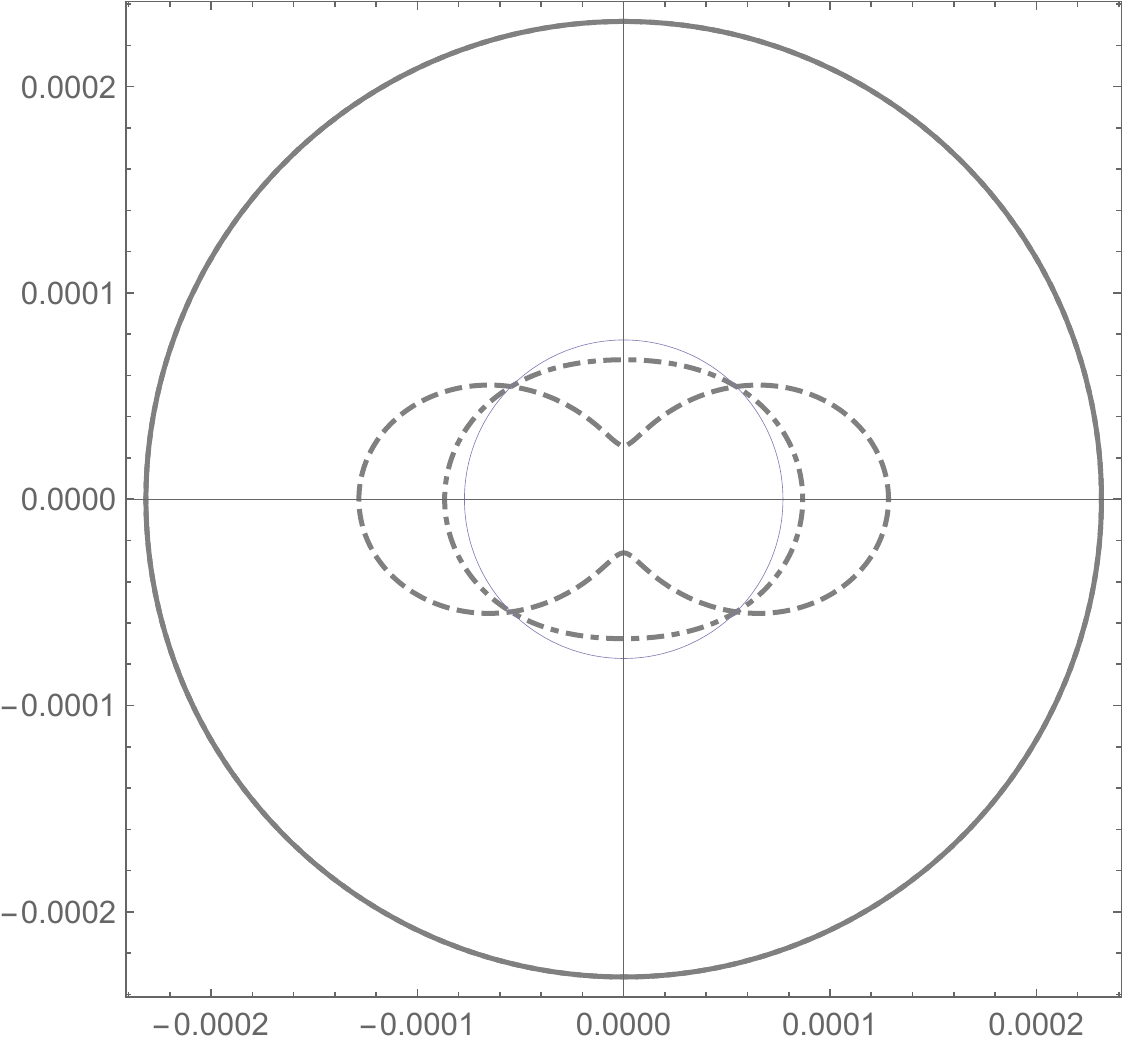}
\caption{Polar diagrams of $\u_{1}(\theta)=\w_{1}(\theta)$, continuous line, $(\v_{1})_1(\theta)$, dashed line, $(\v_{2})_1(\theta)$, dotdashed line for the second laminate; the thin line is the zero-line: inside it, a diagram represents a negative value.}
\label{fig:16}
\end{figure}

\begin{figure}[h]
\center
\includegraphics[width=.6\textwidth]{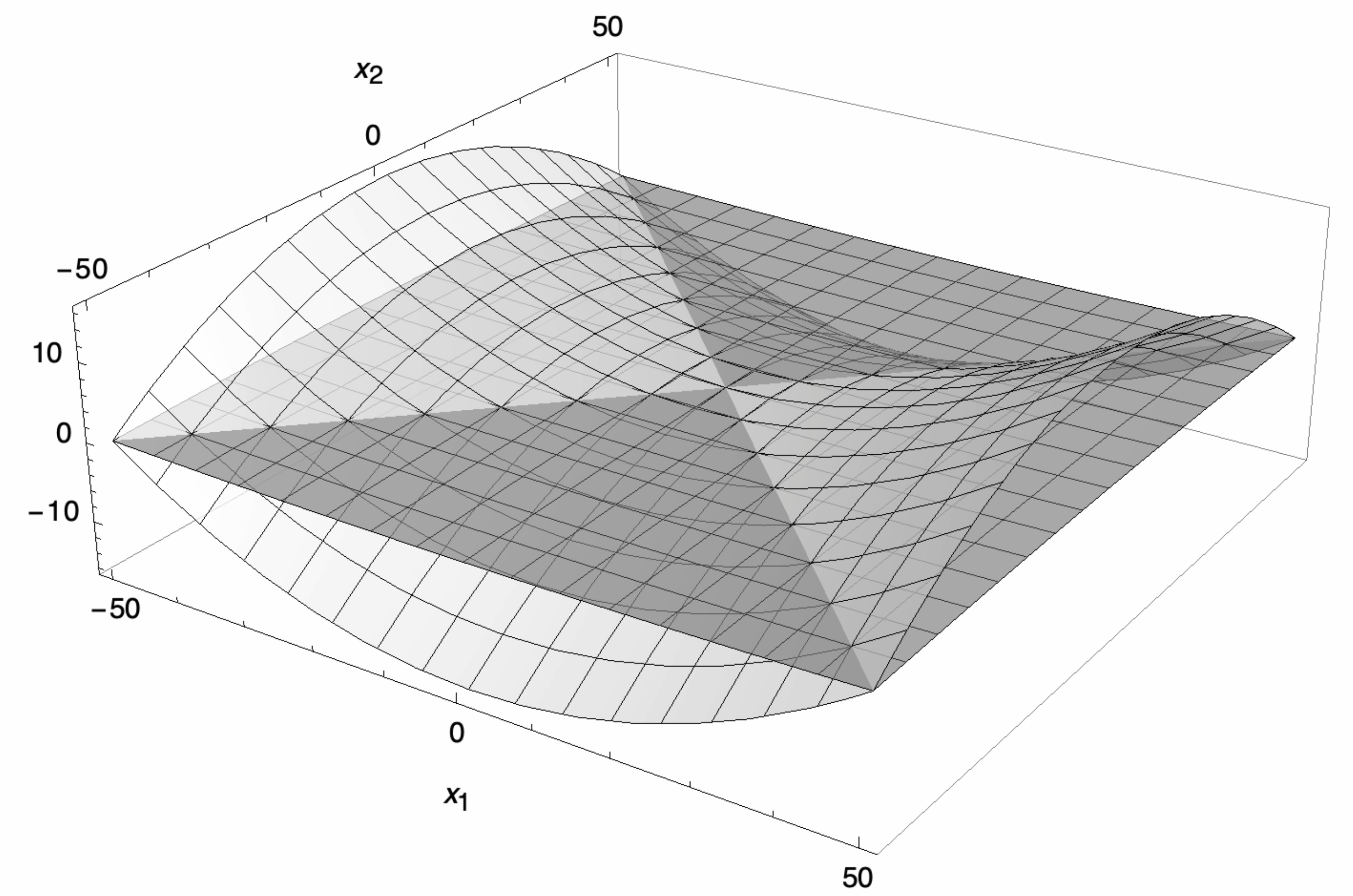}
\caption{Comparison between the two bent shapes of the first (light grey) and second (dark grey) laminate for $t=50^\circ$C, $\nabla t=0$.}
\label{fig:17}
\end{figure}

\begin{figure}[h]
\center
\includegraphics[width=.6\textwidth]{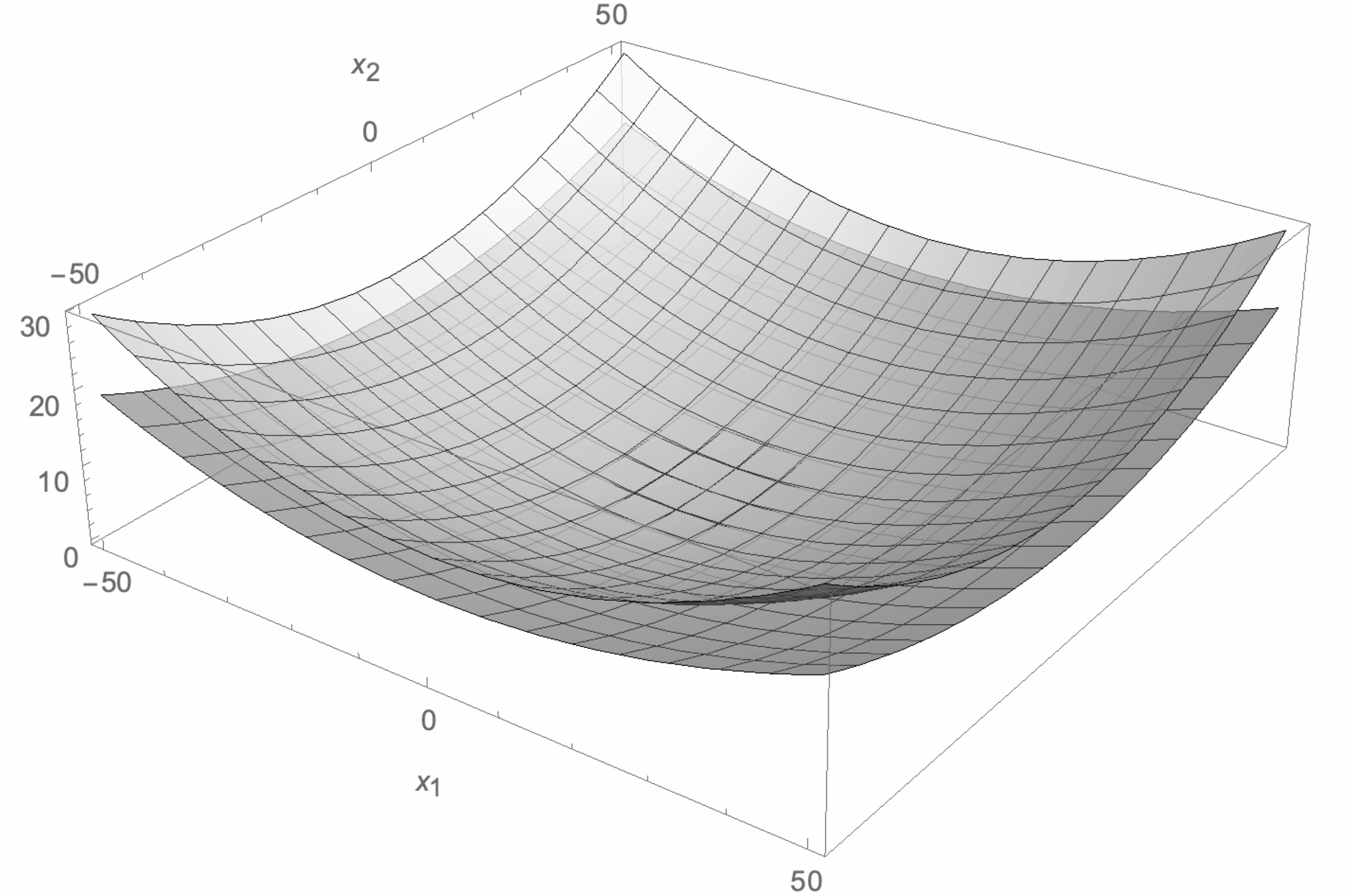}
\caption{Comparison between the two bent shapes of the first (light grey) and second (dark grey) laminate for $t=0,\nabla t=50^\circ$C/mm.}
\label{fig:18}
\end{figure}

Moreover, the curvature tensor $\bk$ is isotropic, because $\w$ is isotropic:
\be
\bk=\frac{(R_0^A+T_0) T_\gamma-6 (-1)^\lambda {R_1^B}^2 \rho}{4\left[(R_0^A+T_0) T_1-6 {R_1^B}^2\right]}\left\{\begin{array}{c}\nabla t\\\nabla t\\0   \end{array}\right\}=\left\{\begin{array}{c}5.21\\5.21\\0   \end{array}\right\}\nabla t\times10^{-4} \mathrm{mm}^{-1}.
\ee
So, as an effect of $\nabla t$, the plate bends like a part of a sphere, Fig. \ref{fig:12}. Unlike $\v_1$, $\w$ does not depend on $h$: the change of number of plies does not affect the bending of the plate.
\\
{\bf Second laminate}: this is the last sequence in eq. (\ref{eq:sequenze}); it is completely asymmetric. As $\A$ and $\U$ do not depend on the stacking sequence, because $a_k=1/n\ \forall k$ in eq. (\ref{eq:coefABCD}), but only on orientations, they are the same of the previous laminate. Then, because the sequence is quasi homogeneous,  $\D=\A$ and $\W=\U$, so  also $\D$ and $\W$ do not change with respect to the previous example. What really changes is $\B$ and $\V$, which makes change also $\Ac,\Bc,\Dc,\u,\v_1,\v_2$ and $\w$. The new values for these tensors are:\medskip
\begin{itemize}
\item $\B=\left[\begin{array}{ccc}-4.76&0&0\\0&4.76&0\\0&0&0\end{array}\right]\times10^3$ MPa, 
$R_1^B=1.19\times10^3$ MPa;\bigskip
\item $\Ac\hspace{-1mm}=\hspace{-1mm}\Dc\hspace{-1mm}=\hspace{-1mm}\left[\begin{array}{ccc}1.05&-0.03&0\\-0.03&1.05&0\\0&0&6.97\end{array}\right]\hspace{-.5mm}\times\hspace{-0.5mm}10^{-5}$ MPa$^{-1}$, $\Bc\hspace{-1mm}=\hspace{-1mm}\left[\begin{array}{ccc}3.47&0&0\\0&-3.47&0\\0&0&0\end{array}\right]\hspace{-0.5mm}\times\hspace{-0.5mm}10^{-5}$ MPa$^{-1}$,\\
$t_0^A=t_0^D=2.01\times10^{-5}$ MPa$^{-1}$, $t_1^A=t_1^D=2.54\times10^{-6}$ MPa$^{-1}$,\\
$r_0^A=r_0^D=1.47\times10^{-5}$ MPa$^{-1}$, $r_1^A=r_1^D=0, \phi_0^A=\phi_0^D=45^\circ,\phi_1^A=\phi_1^D=0$;\bigskip
\item $\V=\left\{\begin{array}{c}0.46\\-0.46\\0   \end{array}\right\}\ ^\circ$C$^{-1}$MPa, 
 $T^V=0,R^V=0.46\ ^\circ$C$^{-1}$MPa; \bigskip
\item $\u=\w=\left\{\begin{array}{c}1.54\\1.54\\0   \end{array}\right\}\times10^{-4}\ ^\circ$C$^{-1}$, $\v_1=\left\{\begin{array}{c}5.11\\-5.11\\0   \end{array}\right\}\times10^{-5}\ ( ^\circ$C mm)$^{-1}$,\\ $\v_2=\left\{\begin{array}{c}9.59\\-9.59\\0   \end{array}\right\}\times10^{-6}\ ^\circ$C$^{-1}$mm,\\
$t^u=t^w=1.54\times10^{-4}\ ^\circ$C$^{-1}$, $r^u=r^w=0, t^v_1=t^v_2=0,$\\
$r^v_1=5.11\times10^{-5}\ ( ^\circ$C mm)$^{-1}$, $r^v_2=9.59\times10^{-6}\ ^\circ$C$^{-1}$mm.\medskip
\end{itemize}
The polar diagrams for these new tensors are shown in Figs. \ref{fig:13} to \ref{fig:16}. For the sake of comparison, Figs. \ref{fig:17} and \ref{fig:18} show the bent shape of the plate for the two cases of $t=50^\circ$C, $\nabla t=0$ and $t=0,\nabla t=50^\circ$C/mm, respectively. 
It is apparent that, globally, the effects of the coupling tensor $\v_1$ decrease: in Fig. \ref{fig:17} the deformation of the second laminate is very small compared to that of the first one. Also the direct effect of tensor $\w$ decreases, Fig. \ref{fig:18}, but not as much as that of the coupling tensor $\v_1$.

\section{Conclusion}
A general analysis of the thermoelastic coupling properties and of their effects for anisotropic laminates has been presented. The use of the polar formalism has allowed to obtain some general results, both for the qualitative and quantitative aspects. It has been proved the existence of a physical constant, the ratio $\rho$, for the basic layer and for all the elastic tensors $\A,\B,\D$ and the corresponding thermoelastic ones, $\U,\V,\W$. Also, the two classes of laminates with $\V=\gr{O}$, allowing for fabricating thermally stable laminates, and those with $\V\neq\gr{O}$, have been deeply examined, and some interesting cases for practical applications analyzed. The examples given show the importance of  thermal coupling effects for what concerns the deformation of a coupled laminate. 

It is worth noting that the whole study and results of this paper can be applied as well to another problem, similar to that of temperature changes: the problem of moisture absorption, \cite{CompStruct07,vannucci_libro}: though this phenomenon happens on a much longer time scale than that of thermal expansion/bending, it is ruled exactly by the same type of physical laws: tensors $\U,\V,\W,\u,\v_1,\v_2$ and $\w$ are defined exactly in the same way and have exactly the same physical meaning, but they concern moisture absorption and not temperature changes.

The use of TQHCL is particularly interesting for finding thermally coupled laminates, but it is not the only one; another possible approach is to use fully orthotropic coupled laminates as described in \cite{CompStruct05}, where the way to find them is also indicated, or to make use of numerical approaches based upon the statement of an adequate optimization problem, like described, e.g., in \cite{vannucci06,vannucci09algo,vincenti10}. Finally,  the design of structures composed by thermally coupled laminates is a rather complicate task; the use of suitable numerical approaches, like the one described in \cite{montemurro12jota1,montemurro12cs}, tackling, on the one hand, the design of the overall structural properties and, on the other hand, the stacking sequence, is necessary, research is going on in this direction.

\bibliographystyle{plain} 
\bibliography{Biblio}
\end{document}